\documentclass[aos]{imsart}

\RequirePackage{amsmath,amsthm,amsfonts,amssymb}
\RequirePackage[authoryear]{natbib}
\RequirePackage[colorlinks,citecolor=blue,urlcolor=blue]{hyperref}
\RequirePackage{comment}
\RequirePackage{booktabs}
\usepackage{subfig}
\usepackage{graphicx}

\startlocaldefs
\numberwithin{equation}{section}
\theoremstyle{plain}

\newtheorem{theorem}{Theorem}[section]

\newtheorem{corollary}[theorem]{Corollary}
\theoremstyle{remark}

\newtheorem{remark}[theorem]{Remark}
\newtheorem{example}{Example}
\newtheorem{assumption}{Assumption}

\newcommand{\C}[1]{\mathcal{#1}}

\newcommand{\gh}{\Hat{\gamma}}
\newcommand{\ah}{\Hat{\alpha}}
\newcommand{\Vh}{\Hat{V}}
\renewcommand{\th}{\Hat{\theta}}
\DeclareMathOperator{\E}{E}

\usepackage[suppress]{color-edits}
\addauthor{vs}{red}

\endlocaldefs

\begin{document}

\begin{frontmatter}
\title{Automatic Debiased Machine Learning \\ Via Riesz Regression}
\runtitle{AutoDML via Riesz Regression}

\begin{aug}
\author[A]{\fnms{Victor}~\snm{Chernozhukov}\ead[label=e1]{vchern@mit.edu}},
\author[A]{\fnms{Whitney K.}~\snm{Newey}\ead[label=e2]{wnewey@mit.edu}}
\author[A]{\fnms{Víctor}~\snm{Quintas-Martínez}\ead[label=e3]{vquintas@mit.edu}}
\and
\author[B]{\fnms{Vasilis}~\snm{Syrgkanis}\ead[label=e4]{vsyrgk@stanford.edu}}
\address[A]{Department of Economics,
MIT\printead[presep={,\ }]{e1,e2,e3}}

\address[B]{Department of  Management Science and Engineering,
Stanford University\printead[presep={,\ }]{e4}}
\end{aug}

\begin{abstract}
A variety of interesting parameters may depend on high dimensional regressions. 
Machine learning can be used to estimate such parameters. 
However estimators based on machine learners can be severely biased by regularization and/or model selection. 
Debiased machine learning uses Neyman orthogonal estimating equations to reduce such biases. 
Debiased machine learning generally requires estimation of unknown Riesz representers. 
A primary innovation of this paper is to provide Riesz regression estimators of Riesz representers that depend on the parameter of interest, rather than explicit formulae, and that can employ any machine learner, including neural nets and random forests.
End-to-end algorithms emerge where the researcher chooses the parameter of interest and the machine learner and the debiasing follows automatically.
Another innovation here is debiased machine learners of parameters depending on generalized regressions, including high-dimensional generalized linear models.
An empirical example of automatic debiased machine learning using neural nets is given. 
We find in Monte Carlo examples that  automatic debiasing sometimes performs better than debiasing via inverse propensity scores and never worse.
Finite sample mean square error bounds for Riesz regression estimators and asymptotic theory are also given. 
\end{abstract}



\end{frontmatter}


\section{Introduction}
Many parameters of interest depend on regressions. 
Examples include treatment effects, regression decompositions, and policy effects. 
Often, a regression may be high dimensional, depending on many variables. 
For example there may be many covariates for treatment effects.   Machine learning methods such as neural nets, random forests, and Lasso can be used to estimate parameters of interest that depend on high dimensional regressions.

A general problem with estimating parameters of interest using machine learning is that machine learners are biased by regularization and/or model selection. 
This bias may pass through when the learner is plugged into a formula for a parameter of interest and make the parameter estimator highly biased. 
This problem can be avoided by using Neyman orthogonal estimating equations where machine learners have zero first-order effect. 
Cross-fitting, a form of sample splitting, can also help.

The orthogonal estimating equations for regressions depend on a Riesz representer  $\alpha_0$ that must be estimated. 
The primary innovation of this paper is to provide an automatic estimator of $\alpha_0$ that uses only the definition of the parameter of interest and the regression but does not require knowing a formula for  $\alpha_0$.  
We give an objective function with expectation that is minimized at $\alpha_0$ that depends only the parameter of interest. 
We refer to minimization of this objective function as a Riesz regression, being equivalent to minimizing the expected squared deviation from  $\alpha_0$. 
Neural nets, random forests, and other methods can be used for this Riesz regression. 
Using the Riesz regression estimator in the bias correction completes an algorithm that 1) specifies the parameter of interest;
2) specifies a learner of the unknown regression; and 3) uses the Riesz regression estimator of $\alpha_0$ determined by steps 1) and 2).
 
A second innovation of this paper is to construct and derive properties of estimators that depend on generalized regressions, which minimize an expected loss over some linear set of functions. 
These generalized regressions include conditional means, least squares projections, functions that minimize quasi-likelihoods, and quantile regressions.
Debiasing for generalized regressions depends on a weighted version of the Riesz representer. 
We give a weighted Riesz regression that only uses the parameter of interest and the generalized regression for bias correction.

A third contribution of this paper is finite sample mean square error bounds for Reisz regressions. 
These bounds are obtained using the critical radius of functions of $\alpha$ on which the objective function depends and approximation error bounds for the unknown $\alpha_0$. 
A fourth contribution is convergence rates for neural net Riesz regressions.
These are based on known results on critical radius and approximation error for neural nets and the finite sample bounds given here.

In work that followed up on the first version of this paper (\cite{chernozhukov2022riesznet}) we found that using Riesz regressions to debias neural net and random forest estimators of the average treatment effect was much more accurate than state of the art methods based on inverse propensity score weighting, in a Monte Carlo study.
Both the automatic neural net and random forest debiasing also led to accurate confidence intervals in those experiments.

Automatic debiasing for Lasso and reproducing kernel Hilbert space regressions was previously given by \cite{chernozhukov2022automatic} and \cite{singh2022kernel} respectively. 
The estimator of $\alpha_0$ given here goes beyond these to provide automatically debiasing for generalized regressions based on neural nets, random forests, and other machine learners. 
These innovations allow researchers to use any of a wide variety of automatically debiased machine learners learners to estimate parameters of interest that depend on generalized regressions. 
For example, automatic debiased machine learning with neural nets could be especially useful for parameters that depend on  high dimensional, nonlinear generalized regressions.

This paper builds on recent work on Neyman orthogonal scores and debiased machine learning. 
We use model free orthogonal estimating equations like those of \cite{chernozhukov2022locally} 
that are the sum of an identifying moment function and a bias adjustment (influence function) term for generalized regression from \cite{ichimura2022influence}.
Those papers did not give the Riesz regression.
Finite sample mean square error bounds for a general learner of $\alpha_0$ are obtained by applying results of \cite{foster2019orthogonal} that characterize error bounds in terms of critical radius and approximation. 
The rate of convergence for neural net Riesz regression use critical radius and approximation rate results given in \cite{farrell2021deep}. 
Additional neural net rate conditions could be obtained using \cite{yarotsky2018optimal}. 
The learner of $\alpha_0$ differs from those of \cite{farrell2021deep, farrell2021deep2} in using the Riesz regression rather than a known form for $\alpha_0$.

We also build upon ideas in classical semi- and nonparametric learning theory with low dimensional regressions using traditional smoothing methods (\cite{van1991differentiable, bickel1993efficient, newey1994asymptotic, robins1995semiparametric, van2000asymptotic}),
that do not apply to  machine learners.  
The orthogonal estimating equations given in \cite{chernozhukov2022locally} and used here build on previous work on nonparametric orthogonal moment functions by  \cite{levit1975efficiency,hasminskii1978nonparametric,bickel1988estimating,newey2004twicing}.
Targeted maximum likelihood (\cite{van2006targeted}) based on machine learners has been considered by \cite{van2011targetedbook}  and large sample theory given by \cite{luedtke2016optimal}.  

In section \ref{sec:linear} we give the Reisz regression and an automatic debiased machine learning algorithm for parameters that depend linearly on a nonparametric regression, including examples. Section \ref{sec:linear} also gives finite sample mean square error bounds for the general and neural net Riesz regression and asymptotic inference results for parameters that are linear functionals of a nonparametric regression. 
Section \ref{sec:multip} extends the estimation methods and theory to nonlinear functionals of generalized regressions. 
In section \ref{sec:empir} we illustrate the usefulness of our methods with an empirical application. 
Section \ref{sec:simul} presents the results of our simulation exercises.

\section{Average Linear Effects for a Conditional Mean}\label{sec:linear}
To highlight the innovation provided by the Riesz regression, we first consider average linear effects that depend on a conditional mean. 
In section \ref{sec:ext} we consider the general setting of nonlinear functions of generalized regressions.
\subsection{Parameters of Interest}
We consider data that consists of i.i.d. observations $W_1, \ldots, \allowbreak W_n$, each having CDF $F_0$. A data observation $W$ includes an outcome variable $Y$ and regressors $X$. In this section, we focus on parameters that depend on the conditional mean of $Y$ given $X$. We will denote a possible such regression function by $\gamma$, with $\gamma_0(x) = \E[Y \mid X = x]$ being the true regression function.

The parameter of interest $\theta_0$ has the form
\begin{equation}\label{eq:linmom}
    \theta_0 = \E[m(W, \gamma_0)],
\end{equation}
where $m(w, \gamma)$ is a functional that depends on a data observation $w$ and a possible regression function $\gamma$. For now, we assume that the expected functional $\gamma \mapsto \E[m(W, \gamma)]$ is linear and continuous in $\gamma$, meaning that there is a constant $C>0$ with $|\E[m(W, \gamma)]|^{2}\leq C\E[\gamma(X)^2]$ for all $\gamma$ with $\E[\gamma(X)^2] < \infty$. Under this assumption, there exists a function $v_m$ with $\E[v_m(X)^2] < \infty$ such that 
\begin{equation}\label{eq:rr}
    \E[m(W, \gamma)] = \E[v_m(X)\gamma(X)] \quad \text{for all } \gamma \text{ with } \E[\gamma(X)^2] < \infty.
\end{equation}
The existence of this $v_m$ follows from the Riesz representation theorem, and it is equivalent to the semiparametric variance bound for $\theta_0$ being finite (see \cite{newey1994asymptotic, hirshberg2021augmented, chernozhukov2019global}). For these reasons $v_m$ is often referred to as the Riesz representer. 
In this Section, where the parameter of interest depends on a nonparametric regression, the Riesz representer  $v_m = \alpha_0$ needs to be estimated for debiased machine learning. 
In section \ref{sec:ext}, where  $\gamma$ may be a generalized regression, $\alpha_0$ will be a weighted version of the Riesz regression.

There are a variety of important, empirically relevant parameters of interest that have this form. 
We illustrate with some familiar examples to help highlight and motivate the Riesz regression: 

\begin{example}[Average Treatment Effect \label{ex:ate}] Suppose that $X=(D,Z)$ where $D$ is a binary treatment indicator, and $Z$ are covariates. The parameter of interest is $\theta_0$ in equation \eqref{eq:linmom} with 
\[m(W,\gamma)=\gamma(1,Z)-\gamma(0,Z).\]
If $Y = D \cdot Y(1) + (1 - D) \cdot Y(0)$, where potential outcomes $(Y(1), Y(0))$ are conditionally independent of treatment $D$ given covariates $Z$, then this object is the Average Treatment Effect or ATE (\cite{rosenbaum1983central}). In this example the Riesz representer is 
\[ \alpha_0(X) =  \frac{D}{\pi_0(Z)} - \frac{1-D}{1-\pi_0(Z)},\] where $\pi_0(z) = \Pr(D = 1 \mid Z = z)$ is the propensity score. Here $\alpha_0(X)$ is the difference of the \cite{horvitz1952generalization} weights for treated and untreated and  $\E[\alpha_0(X)^2] < \infty$ if and only if $\E[\pi_0(Z)^{-1}(1-\pi_0(Z))^{-1}] < \infty$.
\end{example}

\begin{example}[Average Marginal Effect] \label{ex:deriv} Suppose again that $X=(D,Z)$ where $D$ is now a continuous treatment or policy variable, and $Z$ are covariates. The parameter of interest is $\theta_0$ in equation \eqref{eq:linmom} with
\[m(W,\gamma) = \partial_d \gamma(X) ,\]
where we denote $\partial_d g(x) = \partial g(x) / \partial d$ for any function $g$. This object can be interpreted as an average treatment effect for continuous treatment $D$, see \cite{imbens2009identification}. Here the Riesz representer is
\[\alpha_0(X)=-\partial_d \ln f_0(X), \]  where $f_0(X)$ is the (true) joint probability density function (pdf) of $X$. Here equation \eqref{eq:rr} follows by integrating by parts, and then multiplying and dividing by $f_0(x)$. 
\end{example}

\begin{example}[Average Policy Effect] Suppose that $\gamma_0$ does not vary with the distribution of $X$. Then, the average effect of a conterfactual shift in the distribution of regressors, from a known distribution with pdf $g_0$ to another known distribution with pdf $g_1$ is the $\theta_0$ of equation \eqref{eq:linmom} with 
\[m(W,\gamma) = \int \gamma(x) (g_1(x) - g_0(x)) dx = \E\left[\frac{g_1(X) - g_0(X)}{f_0(X)} \gamma(X)\right],\] where $f_0(x)$ is the (true) pdf of $X$ in the data. 
Here the Riesz representer is
\[\alpha_0(X)= \frac{g_1(X) - g_0(X)}{f_0(X)}, \]  
with equation \eqref{eq:rr} following from the second equality in the expression for $m(W,\gamma)$.
\end{example}


\subsection{Estimation \label{sec:estim}}
We will base estimation of $\theta_0$ on a Neyman orthogonal estimating equation, i.e. score, where first step estimation has zero first order effects, that is also doubly robust in having expectation zero if either $\gamma=\gamma_0$ or $\alpha=\alpha_0$. This score is
\begin{equation}\label{eq:psi}
\psi(w, \gamma, \alpha, \theta) = m(w, \gamma)- \theta + \alpha(x)(y - \gamma(x)),
\end{equation}
as in \cite{chernozhukov2022locally}, where taking expectations gives, for any $\alpha, \gamma$,
\begin{align}
\E[\psi(W,\gamma,\alpha,\theta_0)] & = \E[m(W,\gamma)] - \theta_0 + \E[\alpha(X)(Y-\gamma(X))] \notag \\
& = \E[m(W,\gamma-\gamma_0)] - \E[\alpha(X)(\gamma(X) - \gamma_0(X))] \notag
\\ & = -\E[(\alpha(X)-\alpha_0(X))(\gamma(X) - \gamma_0(X))].
\end{align}
Here we see that at the true parameter value $\theta_0$, the expectation of the score $\psi(W,\gamma,\alpha,\theta_0)$ differs from zero only to second order and equals zero if either $\gamma=\gamma_0$ or $\alpha=\alpha_0$. Thus the score is Neyman orthogonal and doubly robust in that it has zero expectation if either $\gamma=\gamma_0$ or $\alpha=\alpha_0$.

Estimation of $\alpha_0$ is essential to construction of a debiased machine learner of the parameter of interest. 
The primary innovation of this paper is to give an extremum characterization of $\alpha_0$ and use this to estimate $\alpha_0$. 
This extremum characterization is given by
\begin{align}
\alpha_0 & = \arg\min_\alpha \E[(\alpha_0(X) - \alpha(X))^2] \notag \\
& = \arg\min_\alpha \E[\alpha_0(X)^2 - 2\alpha_0(X)\alpha(X) + \alpha(X)^2] \notag \\
& = \arg\min_\alpha \{ -2\E[v_m(X)\alpha(X)] + \E[\alpha(X)^2] \} \notag \\
& = \arg\min_\alpha \E[-2m(W, \alpha) + \alpha(X)^2], \label{eq:extrem}
\end{align}
where the third equality holds because $\E[\alpha_0(X)^2]$ does not depend on $\alpha$ and $\alpha_0=v_m$, and the fourth equality follows from equation \eqref{eq:rr} with $\gamma=\alpha$. This characterization can be used to estimate $\alpha_0$ by replacing the expectation with the sample average and minimizing over some set of possible $\alpha$ functions.

We call a resulting  estimator of $\alpha_0$ a Riesz regression, motivated by minimization of the least squares objective function in equation \eqref{eq:extrem}. 
This estimator is \textit{automatic} in dependng only on the function $m(w,\gamma)$ that determines the parameter of interest and in not requiring the form of $\alpha_{0}$. 
In particular, this method does not depend on plugging in non-parametric estimates of components of $\alpha_{0}$. 
This feature is useful when $\alpha_0$ does not have a simple form. 
For causal parameters such as those of Examples 1-3, the Riesz regression avoids inverting a learner of a conditional probability or a pdf.  
Instead, the Riesz regression learns $\alpha_{0}$ directly. 

Our estimation strategy for the parameter of interest is to combine the Reisz regression estimator of $\alpha_0$ and estimation of $\gamma_0$ in the Neyman orthogonal score with the use of cross-fitting to reduce overfitting bias.\footnote{See \cite{newey2018cross} for more on the advantages of cross-fitting.} 
The outline of our estimation strategy is as follows:

\begin{enumerate}
    \item Partition the set of data indices ${1, \ldots, n}$ into $L$ disjoint subsets of about equal size $I_\ell$, $\ell = 1, \ldots, L$;
    \item For each data fold $\ell = 1, \ldots, L$:
    \begin{enumerate}
        \item Estimate $\gh_\ell \in \mathcal{G}_n$ as a non-parametric regression of $Y$ on $X$ over some class of functions $\mathcal{G}_n$ using observations \textit{not} in $I_\ell$.
        \item Estimate the debiasing function $\ah_\ell$ using observations \textit{not} in $I_\ell$ by minimizing a sample version of the objective function in equation \eqref{eq:extrem} over a set of functions, as in
\[\ah_\ell = \arg\min_{\alpha \in \mathcal{A}_n} 
\bigg[
\sum_{i \notin I_\ell}
\left\lbrace
- 2m(W_i, \alpha) +
\alpha(X_i)^2
\right\rbrace
+ \Lambda_r(\alpha)
\bigg]
\]
Where $\Lambda_r(\alpha)$ is a penalty term and $r$ is a scalar determining the magnitude of penalization.
    \end{enumerate}
    \item Estimate the parameter of interest using the cross-fitted regression and debiasing function in the moment function of equation \eqref{eq:psi} to obtain
\[\th = \frac{1}{n}\sum_{\ell = 1}^L \sum_{i \in I_\ell} \left\lbrace m(W_i, \gh_\ell) + \ah_\ell(X_i) (Y_i - \gh_\ell(X_i))\right\rbrace\]
    \item Estimate the standard error of $\th$ as $\sqrt{\hat{V}/ n}$, where:
\[\Vh = \frac{1}{n}\sum_{\ell = 1}^L \sum_{i \in I_\ell} \left\lbrace m(W_i, \gh_\ell) + \ah_\ell(X_i) (Y_i - \gh_\ell(X_i)) - \th\right\rbrace^2 \]
\end{enumerate}

Our estimation strategy is very general, allowing for any choice of learner $\hat{\gamma}_{\ell}$ and any Riesz regression $\hat{\alpha}_{\ell}$  encoded in the class of functions $\mathcal{A}_n$. 
Special kinds of Riesz regressions have been given in previous literature. 
These include  linear combinations of a dictionary of functions $(b_1(x), \ldots, b_p(x))'$, and $p$ large, with an $L_1$ penalty in the loss function (\cite{chernozhukov2022automatic}), or functions embedded in a reproducing kernel Hilbert space (\cite{singh2021debiased}). \cite{chernozhukov2020adversarial} allowed for the estimation of $\alpha_0$ in arbitrary function spaces, but proposed a computationally harder minimax loss formulation. 
A primary innovation of this paper is to provide the Riesz regression for automatic estimation of $\alpha_0$ and corresponding asymptotic theory. 

As an example, below we will give primitive conditions for a neural net Riesz regression. 
A general neural net takes the form
\[x\overset{f_{1}}{\longmapsto}H^{(1)}\overset{f_{2}}{\longmapsto}\cdots\overset{f_{m}}{\longmapsto}H^{(m)}\]
where $H^{(l)}=\{H_{k}^{(l)}\}_{k=1}^{K_{l}}$ are called neurons, $x$ is the original finite-dimensional input, and the function $f_{l}$ maps one layer of neurons to the next as in%
\[f_{l}:v\overset{f_{l}}{\longmapsto}\{H_{k}^{(l)}(v)\}_{k=1}^{K_{l}%
}:=(1,\{\sigma(v^{\prime}\beta_{k,l})\}_{k=2}^{K_{l}}),\]
where each $\beta_{k,l}$ is a $K_{l-1}$ vector of parameters and $\sigma(u)$ is a nonlinear activation function. We will focus on the case where $\sigma(u)$ is the RELU function $\sigma(u)=\max\{0,u\}$. An important special case is a multilayer perceptron (MLP) network where the number of neurons $K_{l}=K$ is the same for each layer, for which results were recently given by \cite{farrell2021deep}. Sparse versions of this specification, where many of the elements of the coefficient vectors $\beta_{k,l}$ may be zero, have also been considered recently by \cite{schmidt2020nonparametric}. \cite{yarotsky2018optimal} gave other neural net specifications with good approximation properties.

In the setting of Example \ref{ex:ate}, a neural net Riesz regression would be constructed as $\hat\alpha_\ell(d, z) = \alpha(d, z; \hat\beta_\ell)$, where:
\[\hat\beta_\ell = \arg\min_{\beta} 
\bigg[
\sum_{i \notin I_\ell}
\left\lbrace
- 2[\alpha(1, Z_i; \beta) - \alpha(0, Z_i; \beta) ]  +
\alpha(D_i, Z_i; \beta)^2
\right\rbrace
+ \Lambda_r(\beta)
\bigg]
\]
for some penalty function $\Lambda_r(\beta)$ (e.g., L1, L2, or the elastic net). Because $D$ is binary, a convenient neural net architecture in this case could be a bi-headed MLP, $\alpha(d, z; \beta, \delta_0, \delta_1) = d g(z; \beta)'\delta_1 + (1 - d) g(z; \beta)'\delta_0$, where $g(z; \beta)$ is an MLP. An even more flexible specification would be to have $\alpha(d, z; \beta_0, \beta_1) = d g(z, \beta_0) + (1 - d) g(z, \beta_1)$, i.e. an MLP for the case $d = 1$ and another MLP for the case $d = 0$.

For Example \ref{ex:deriv}, a neural net Riesz regression is  $\hat\alpha_\ell(d, z) = \alpha(d, z; \hat\beta_\ell)$, where:
\[\hat\beta_\ell = \arg\min_{\beta} 
\bigg[
\sum_{i \notin I_\ell}
\left\lbrace
- 2[\partial_d \alpha(D_i, Z_i; \beta)]  +
\alpha(D_i, Z_i; \beta)^2
\right\rbrace
+ \Lambda_r(\beta)
\bigg].
\]
In particular, notice that the loss function involves taking a derivative of the neural net with respect to one of the inputs. A convenient parametrization of the neural net in this case is a locally linear function $\alpha(d, z; \phi, \beta) = \phi(d, z)'g(z; \beta)$, where $\phi(d, z)$ is a dictionary of known, differentiable basis functions, and $g(z; \beta)$ is a neural net (e.g. an MLP). In that case, $\partial_d \alpha(d, z; \beta) = [\partial_d \phi(d, x)]'g(z; \beta)$. This approach was used in \cite{chernozhukov2022riesznet} to construct a random forest estimator of $\alpha_0$, exhibiting good performance in Monte Carlo simulations.

\subsection{Large Sample Inference for Linear Effects of Regression}

In this Section, we give mean square convergence rates for learners $\hat{\alpha}_{\ell}$ and $\sqrt{n}$-consistency and asymptotic normality results for the learner $\hat{\theta}$ of the object of interest and its asymptotic variance estimator $\hat{V}$. We first derive convergence rates for $\hat{\alpha}_{\ell}.$

\subsubsection{Convergence Rates for $\hat{\alpha}_{\ell}$ \label{sec:largesampalpha}}
In this subsection we suppress the $\ell$ subscript for notational convenience.
We consider the problem of estimating
\[
\alpha_{0}=\arg\min_{\alpha}\E[-2m(W,\alpha)+\alpha(X)^{2}],
\]
where we have used the extremum characterization of $\alpha_0$ in equation \eqref{eq:extrem}. For any random variable $a(W)$ let $\left\Vert a\right\Vert
=\sqrt{\E[a(W)^{2}]}$ and $\left\Vert a\right\Vert _{\infty}=\sup
_{w\in\mathcal{W}}\left\vert a(w)\right\vert.$  \vsedit{For simplicity of exposition we will only consider the case where the estimate is defined over a growing sieve space $\mathcal{A}_n$ and no regularization is used, i.e. $\Lambda_r(\alpha)=0$:
\begin{equation}
\hat{\alpha}=\arg\min_{\alpha\in\mathcal{A}_{n}}\sum_{i=1}^n \left\{  -2m(W_{i},\alpha) + \alpha(X_{i})^{2}\right\},%
\end{equation}
Our estimation rate can easily be extended to regularized estimation with appropriate regularization weight.}
We assume that $m(W,\alpha)$ is mean square continuous in the following sense:
\begin{assumption}\label{ass1}
For some $M > 0$ it is the case that $\E[m(W,\alpha)^{2}] \leq M\left\Vert \alpha\right\Vert ^{2}$. 
\end{assumption}

Define:
\begin{align*}
{\mathrm{star}}(\mathcal{A}_n-\alpha_{0})=~  &  \{x\rightarrow
\xi\,(\alpha(x)-\alpha_{0}(x)):\alpha%
\in\mathcal{A}_n,\text{ }\xi\in\lbrack0,1]\}\\
{\mathrm{star}}(m\circ\mathcal{A}_n-m\circ\alpha_{0})=~  &
\{w\rightarrow\xi\,(m(W,\alpha)-m(W,\alpha_{0}%
)):\alpha\in\mathcal{A}_n,\text{ }\xi\in\lbrack0,1]\}
\end{align*}

\begin{assumption}\label{ass2}
$\left\Vert f\right\Vert _{\infty}\leq1$ for all $f\in{\mathrm{star}}(\mathcal{A}_n-\alpha_{0})$ and $f\in {\mathrm{star}}(m\circ\mathcal{A}_n-m\circ\alpha_{0})$.
\end{assumption}

Define:
\[
\alpha^{\ast}=\arg\min_{\alpha\in\mathcal{A}_n}\E[-2m(W,\alpha)+\alpha(X)^{2}]
\]
to be the best approximation of $\alpha_{0}$ by an element of $\mathcal{A}_n$. 

\begin{theorem}\label{thm1}
Let $\delta_{n}$ be an upper bound on the critical radius of $\mathrm{star}(\mathcal{A}_n-\alpha_{0})$ and $\mathrm{star}(m\circ \mathcal{A}_n-m\circ\alpha_{0})$. If Assumptions \ref{ass1} and \ref{ass2} are
satisfied then it follows that with probability $1-\zeta$, for some universal constant $C$,
\[
\left\Vert \hat{\alpha}-\alpha_{0}\right\Vert ^{2}\leq C\left(M\delta_{n}%
^{2}+\left\Vert \alpha^{\ast}-\alpha_{0}\right\Vert ^{2}+\frac{M\ln(1/\zeta
)}{n}\right).
\]
\end{theorem}
See e.g. \cite{foster2019orthogonal} for the definition of the critical radius used in the statement this result.
To use Theorem \ref{thm1} to obtain a mean square convergence rate for $\hat{\alpha}$ it is important to know the critical radius and the rate at which $\Vert\alpha_{\ast}-\alpha_{0}\Vert$ shrinks as the approximating set $\mathcal{A}_n$ becomes richer. For example, \cite{farrell2021deep} have recently obtained such results for deep, ReLU neural nets. We can apply their results to obtain a mean square rate for such a learner of $\alpha_{0}$ when $x$ is a $d$ dimensional input for the multilayer perceptron (MLP) network with $m$ layers and width $K.$

The convergence rate depends on the smoothness of the function $\alpha_{0}(x),$ as specified in the following result. Specifically we assume that the support of $X$ is contained in a Cartesian product $\mathcal{X}$ of compact intervals and $\alpha_{0}(X)$ can be extended to a function that is continuously differentiable on $\mathcal{X}$ and has $\beta$ continuous derivatives.

\begin{corollary}\label{cor1}
If (i) the support of $X$ is contained in a Cartesian product of compact intervals and $\alpha_{0}(X)$ can be extended to a function that is continuously differentiable with $\beta$ continuous derivatives; (ii) $\mathcal{A}_n$ is an MLP network with $d$ inputs, width $K$, and depth $m$ with $K\to\infty$ and $m\to\infty$; (iii) $m\circ\mathcal{A}_n$ is representable as such a network; then there is $C>0$ such that, for any $\varepsilon >0$, 
\[
\left\Vert \hat{\alpha}-\alpha_{0}\right\Vert^{2}=O_{p}(K^{2}m^{2}\ln
(K^{2}m)\ln(n)/n+[Km\sqrt{\ln(K^{2}m)}]^{-2(\beta/d)+\varepsilon}).
\]
\end{corollary}
When $\alpha_{0}$ is smooth enough, in that $\beta$ is large enough, the upper bound on $\left\Vert \hat{\alpha} -\alpha_{0}\right\Vert $ in Corollary \ref{cor1} gives a mean square convergence rate that can be close to, but less than $n^{-1/2}$. Such rate can be obtained by choosing the width $K$ and depth $m$ to approximately balance the two terms in Corollary \ref{cor1}, $K \asymp n^{\frac{d}{2(\beta + d)}}\ln^2(n)$, $m \asymp \ln(n)$, as in \cite{farrell2021deep}, in which case
\[
\left\Vert \hat{\alpha}-\alpha_{0}\right\Vert^{2} = O_p\left(n^{-\frac{\beta}{\beta + d} }\ln^8(n)\right).
\]
Faster rates could be obtained using the neural nets of \cite{yarotsky2018optimal} or the sparse neural nets of \cite{schmidt2020nonparametric}. We focus on Corollary \ref{cor1} for the MLP neural net because it is a widely used architecture in practice, and because the rates obtained are fast enough for the estimators of the parameter of interest to be asymptotically normal.

\subsubsection{Large Sample Inference for $\theta_{0}$}
We use additional regularity conditions to show asymptotic normality of $\hat{\theta}$ and consistency of the asymptotic variance estimator $\hat{V}$. We will first give a general result for $\hat{\theta}$ that applies to any $\hat{\alpha}_{\ell}$ and does not rely on Theorem \ref{thm1} for a convergence rate for $\hat{\alpha}_{\ell}$. Similarly, any regression learner $\hat{\gamma}_\ell$ can be used here as long as its mean-square convergence rate is fast enough, as formalized below. Such convergence rate results are available for shallow (\cite{chen1999improved}) and deep (\cite{yarotsky2018optimal, schmidt2020nonparametric, farrell2021deep}) neural nets, random forests (\cite{syrgkanis2020estimation}), LASSO (\cite{bickel2009simultaneous}), boosting (\cite{luo2022highdimensional}) and other high-dimensional methods.

The following assumption imposes a few additional regularity conditions. Let $\sigma_0^2(X) = \E[(Y - \gamma_0(X))^2 \mid X]$ denote the conditional variance of $Y$ given $X$.

\begin{assumption}\label{ass3} 
$\alpha_0(X)$ and $\sigma^2_0(X)$ are bounded and $\E[m(W, \gamma_0)^2] < \infty$.
\end{assumption}

Next, we require mean square consistency of $\hat{\gamma}_{\ell}$ and $\hat{\alpha}_{\ell},$ that the product of their mean-square convergence rates is smaller than $n^{-1/2},$ and a boundedness condition for $\hat{\alpha}_{\ell}$.

\begin{assumption}\label{ass4}
(i) $\left\Vert \hat{\gamma}_{\ell}-\gamma_{0}\right\Vert \overset{p}{\to}0$ and $\left\Vert\hat{\alpha}_{\ell}-\alpha_{0}\right\Vert \overset{p}{\to}0$; (ii) $\sqrt{n}\left\Vert \hat{\gamma}_{\ell}-\gamma_{0}\right\Vert \left\Vert \hat{\alpha}_{\ell}-\alpha_{0}\right\Vert\overset{p}{\to}0$; (iii)  $\hat{\alpha}_{\ell}(X)$ is bounded.
\end{assumption}

Part (i) implies that both $\hat{\gamma}_\ell$ and $\hat{\alpha}_\ell$ are consistent in mean square. Part (ii) captures an important tradeoff between the rates of convergence for $\hat{\gamma}_\ell$ and $\hat{\alpha}_\ell$. In settings where the regression can be estimated at a relatively fast rate of convergence, the learner for the debiasing function can converge more slowly, and vice versa, as long as the product of their mean-square convergence rates vanishes faster than $n^{-1/2}$. The results we have obtained for the neural net learner $\hat{\alpha}_{\ell}$ can be used to verify these conditions and we do so in Corollary \ref{cor2} to follow. The mean square convergence of $\hat{\gamma}_{\ell}$ is a primitive condition for this paper and allows use of a wide variety of $\hat{\gamma}_{\ell}$ in the construction of the estimator.

We have the following large sample inference result under these conditions.
\begin{theorem}\label{thm2}
    If Assumptions \ref{ass1}, \ref{ass3} and \ref{ass4} are satisfied, then
    \[
    \sqrt{n}(\hat{\theta}-\theta_{0})\overset{d}{\to}N(0,V)\quad\text{and}\quad\hat{V}\overset{p}{\to}V,
    \]
where $\hat V$ is the variance estimator defined in subsection \ref{sec:estim} and $V=\E[\{m(W,\gamma_{0})-\theta_{0}+\alpha_{0}(X)(Y - \gamma_0(X))\}^{2}]$.
\end{theorem}

Next we use Theorem \ref{thm1} and Corollary \ref{cor1} to formulate regularity conditions when
$\hat{\alpha}_{\ell}$ is the neural net learner of $\alpha_{0}$ in section \ref{sec:estim}. Let
\[
\epsilon_{\alpha n}^2=K^{2}m^{2}\ln(K^{2}m)\ln(n)/n+[Km\sqrt{\ln(K^{2}%
m)}]^{-2(\beta/d)+\varepsilon}.
\]
This $\epsilon_{\alpha n}^2$ is taken from the upper bound for $\left\Vert\hat{\alpha}_{\ell}-\alpha_{0}\right\Vert ^{2}$ in Corollary \ref{cor1} and so characterizes the mean square convergence rate of the automatic neural net learner $\hat{\alpha}_{\ell}$.

\begin{corollary}\label{cor2}
    Suppose that Assumptions \ref{ass1}, \ref{ass2}, \ref{ass3} and the hypotheses of Corollary \ref{cor1} hold. Moreover, suppose that $\left\Vert \hat{\gamma}_{\ell}-\gamma_{0}\right\Vert \overset{p}{\to}0$ , $\epsilon_{an} \to 0$ and $\sqrt{n}\left\Vert \hat{\gamma}_{\ell} - \gamma_{0}\right\Vert \epsilon_{an} \overset{p}{\to}0$. Then, for the neural net learner $\hat{\alpha}_\ell$ of Corollary \ref{cor1}, we have
    \[
    \sqrt{n}(\hat{\theta}-\theta_{0})\overset{p}{\to}N(0,V)\quad\text{and}\quad\hat{V}\overset{p}{\to}V.
    \]
\end{corollary}

\section{Average Effects for Generalized Regressions} \label{sec:ext}
\subsection{Linear  Effects}\label{sec:gener}

In this section we extend the results to parameters that depend on
functions $\gamma_{0}$ other than the conditional mean, that we refer to as
\textit{generalized regressions}. Suppose that $\gamma_0$ is defined as the solution to a general $M$-estimation problem:
\begin{equation}
  \gamma_0 := \arg\min_{\gamma \in \Gamma} \E[\ell(W, \gamma)], \label{eq:M-est}  
\end{equation}
where $\Gamma$ is a closed (in mean square) linear subspace of $L^2(X)$. For example, when $\ell(W, \gamma) = \frac{1}{2} (Y - \gamma(X))^2$ is the square loss and $\Gamma = L^2(X)$, then $\gamma_0(X) = \E[Y \mid X]$ and we recover the case of regression. 

By the first order condition of the minimization problem \eqref{eq:M-est}, $\gamma_0$ satisfies
\begin{equation}
\label{eq:orthotresid}
\E[\rho(W,\gamma_{0})b(X)]=0 \quad \text{for all } b\in\Gamma.
\end{equation}
for some functional $\rho(W, \gamma)$, typically a generalized notion of the (negative) derivative of the loss function $\ell(W, \gamma)$. In the case of regression, we can take $\rho(W, \gamma) = Y - \gamma(X)$ to be the non-parametric residual. For other statistical problems, we will refer to the function $\rho(W, \gamma)$ as a \textit{generalized residual}. The results of this section will apply to any $\gamma_0$ that is identified by an orthogonality condition as in \eqref{eq:orthotresid}, even beyond $M$-estimation problems.

This setting covers many interesting features of the conditional distribution of $Y$ given $X$. First, suppose that $\Gamma = L^2(X)$, so that the functional form of $\gamma_0$ is unrestricted. For example, when $\rho(W ,\gamma) = \tau - 1(Y < \gamma(X))$ for $0 < \tau < 1$, then $\gamma_0(x)$ is the $\tau$-th conditional quantile of $Y$ given $X = x$. When $\rho(W,\gamma)=\lambda(\gamma(X))[Y-\mu(\gamma(X))]$ for a link function $\mu(a)$ and another function $\lambda(a)$, this corresponds to the first order conditions of a generalized linear model \citep{nelder1972generalized}. For binary $Y \in \{0, 1\}$, $\mu(a)$ the standard logistic CDF, and $\lambda(a) \equiv 1$, for instance, this set up corresponds to a (non-parametric) logistic regression, where $\gamma_0(X) = \mu^{-1}(\Pr(Y = 1 \mid X)) = \ln (\Pr(Y = 1 \mid X)/ \Pr(Y = 0 \mid X))$ corresponds to the log-odds.

The set $\Gamma$ could be used to encode parametric or semi-parametric restrictions on $\gamma_0$. One example is $X=(X_{1},X_{2},...)$ and $\Gamma$ the mean square closure of finite linear combinations of $X$. This corresponds to a high (infinite) dimensional, approximately sparse $\Gamma$, where the orthogonality condition is equivalent to $\E[X_{j} \rho(W,\gamma_{0})]=0$ for all $j$. In a case considered also by \cite{hirshberg2021augmented} and \cite{farrell2021deep}, $\Gamma$ is the mean square closure of $\{a(X_{1})+X_{2}^{\prime}b(X_{1})\}$ where $a(X_{1})$ is a scalar function and $b(X_{1})$ a vector of functions, each having unrestricted functional form. We could also take $\Gamma$ to be the mean-square closure of additive functions $a_{1}(X_{1})+a_{2}(X_{2})$, where $X_{1}\,\ $and $X_{2}$ are distinct components of $X$. In some cases, the resulting $\gamma_0$ will have a projection interpretation: for instance, when $\ell(W, \gamma)$ is the square loss, then $\gamma_0 = \arg\min_{\gamma \in \Gamma} \E[(\E[Y \mid X] - \gamma(X))^2]$ is the best approximation to $\E[Y \mid X]$ in $\Gamma$ in the mean-square sense.

For now, we also continue to assume that the parameter of interest has the form $\theta_{0} = \E[m(W,\gamma_{0})]$, where $\gamma\mapsto m(W,\gamma)$ is linear and $\E[m(W,\gamma)]$ is mean square continuous on $\Gamma$. We will relax the linearity assumption in section \ref{sec:multip}. We will extend the results of the previous section by modifying \eqref{eq:psi}. Define, for any $\gamma, \alpha \in \Gamma$, a score:
\begin{equation}
\psi(w,\gamma,\alpha,\theta)=m(w,\gamma)-\theta+\alpha(x)\rho(w,\gamma), \label{eq:orthmomgenreg}
\end{equation}
where we have replaced $y - \gamma(x)$ with the generalized residual $\rho(W,\gamma)$. 

This score satisfies $\E[\psi(W, \gamma_0, \alpha, \theta_0)] = \E[\alpha(X) \rho(W, \gamma_0)] = 0$ for any $\alpha \in \Gamma$ by \eqref{eq:orthotresid}, and hence it is Neyman-orthogonal with respect to $\alpha$. To get Neyman orthogonality with respect to $\gamma$ we need to find a function $\alpha_0$ that satisfies 
\begin{equation}
    \left. \frac{\partial}{\partial r} \E[\psi(W, \gamma_0 + r\delta, \alpha_0)] \right|_{r = 0} = \E[\{v_m(X) + \alpha_0(X)\,v_{\rho}(X)\}\delta(X)] = 0 \quad \forall \delta \in \Gamma, \label{eq:orthogeneralized}
\end{equation}
where $v_m(X)$ is the Riesz representer in \eqref{eq:rr} and, for a scalar $a$,
\begin{align*}
v_\rho(X) := \left. \frac{\partial}{\partial a} \E[ \rho(W, \gamma_0 + a) \mid X] \right|_{a=0},
\end{align*}
that we assume exists. 
We further assume that we can normalize the sign of $\rho(W,\gamma)$ so that $v_{\rho}(X)\leq0$, as will
hold when $\E[\rho(W, \gamma_{0} + a) \mid X]$ is monotonically decreasing in $a$.  For example, when $\rho(W,\gamma)=Y-\gamma(X)$ as in Section 2 we have $v_{\rho}(X)=-1$. Also, when $\rho(W,\gamma)=p-1(Y<\gamma(X))$ then $v_{\rho}(X)=-f_{Y\mid X}%
(\gamma_{0}(X)\mid X),$ the negative of the conditional pdf of $Y$ given $X$
evaluated at $y=\gamma_{0}(X)$.\footnote{Note that since $\gamma_0$ corresponds to the $p$ quantile of $Y\mid X$, if we denote with $\gamma_0(p, X)$ the $p$-th conditional quantile and with $\partial_p\gamma_0(p, X)$ its derivative with respect to $p$, then we have $v_\rho(X) = - (\partial_p\gamma_0(p, X))^{-1}$.}
The Neyman-orthogonality condition above includes $v_\rho(X)$, which was previously equal to $-1$.
Here $v_\rho(X)$ is needed to account for the effect of $\gamma$ on the residual  $\rho(W,\gamma)$.

\begin{remark} 
The orthogonal score will also be doubly robust, in the sense that $\E[\psi(W,\theta_{0},\gamma,\alpha_{0})]=0$ for all $\gamma\in\Gamma,$ if and only if $\E[\alpha_{0}(X)\rho(W,\gamma)]$ is affine in $\gamma$. 
This follows from $\E[m(W,\gamma)]$ being linear in $\gamma$ and from \cite{chernozhukov2022locally}. 
There are many interesting cases where double robustness does not hold, such as conditional quantiles or generalized linear models. 
Even if the score in equation \eqref{eq:orthmomgenreg} is not doubly robust, it will still be orthogonal, enabling $\sqrt{n}$-consistent estimation and asymptotically normal inference on $\theta_{0}$ when $\gamma_{0}$ and $\alpha_{0}$ are estimated by machine learning.
\end{remark}

A key innovation of our work is to note that Equation \eqref{eq:orthogeneralized} can be viewed as the first order condition to the following optimization problem:
\begin{align}
\alpha_0 & = \arg\min_{\alpha\in\Gamma} \E[-2v_{m}(X)\alpha(X)-v_{\rho
}(X)\alpha(X)^{2}]\nonumber\\
& =\arg\min_{\alpha\in\Gamma}\left\{-2\E[v_{m}(X)\alpha(X)]-\E[v_{\rho}(X)\alpha(X)^{2}]\right\}\nonumber\\
& =\arg\min_{\alpha\in\Gamma} \E[-2m(W,\alpha)-v_{\rho}(X)\alpha
(X)^{2}],
\label{eq:wtrieszreg}
\end{align}
where the second equality follows by linearity of expectations, and the third equality follows by  equation \eqref{eq:rr}. 
Thus, $\alpha_{0}$ minimizes the expectation of an objective function that depends on $\alpha$ only through the functional of interest $m(W,\alpha)$ and $\alpha(X)$. 
As with equation \eqref{eq:extrem}, minimizing this objective function does not require any knowledge of the form of $\alpha_0$.

When $v_{\rho}(X)\neq-1$ the $\alpha_0$ will not be the Riesz representer $v_m(X)$.
Instead, $\alpha_0$ can be interpreted as minimizing weighted least squares criterion that depends on the Riesz representer. 
As shown in \cite{ichimura2022influence}, 
\begin{align}
\alpha_{0}  & =\arg\min_{\alpha\in\Gamma} \E\left[-v_{\rho}(X)\left(-\frac{v_{m}(X)}{v_{\rho}(X)} - \alpha(X)\right)^2\right].
\label{eq:wtleastsq}
\end{align}
Thus $\alpha_0(X)$ minimizes a weighted least square criterion with weight $-v_{\rho}(X)$ and the variable being predicted given by $-v_m(X)/v_{\rho}(X)$. 
For this reason we refer to $\alpha_0(X)$ as a weighted Riesz regression. 

Though the objective functions of equations \eqref{eq:wtrieszreg} and \eqref{eq:wtleastsq} differ only by a constant, only equation \eqref{eq:wtrieszreg} possesses the desirable properties that we set out to accomplish of depending solely on known functions of $\alpha$.  
The objective in equation \eqref{eq:wtrieszreg} was not given in \cite{ichimura2022influence}.

\vsedit{In some cases there will be a function $\bar{v}_{\rho}(W)$ such that $\E[\bar{v}_{\rho
}(W)\mid X]=v_{\rho}(X).$ By iterated expectations, the objective function is not
affected by replacing $v_{\rho}(X)$ with $\bar{v}_{\rho}(W),$ because%
\[
\E[-2m(W,\alpha)-v_{\rho}(X)\alpha(X)^{2}]=\E[-2m(W,\alpha)-\bar{v}_{\rho}%
(W)\alpha(X)^{2}].
\]
In practice, it may be easier to minimize the objective function that depends
on $\bar{v}_{\rho}(W)$ to avoid having to estimate $v_{\rho}(X)=\E[\bar{v}_{\rho}(W)\mid X]$.
For this reason, we focus on a sample objective function that depends on an estimator $\hat{v}_\rho(W)$ of $\bar{v}_\rho(W)$ that is allowed to take $W$ as input, instead of just $X$.}

To obtain an estimate of $\alpha_{0}$, we replace the sample criterion in step
(2b) of the algorithm in subsection \ref{sec:estim} with:
\begin{equation}
\hat{\alpha}_{\ell}=\arg\min_{\alpha\in\mathcal{A}_{n}}\sum_{i\notin I_{\ell}%
}\left\{  -2m(W_{i},\alpha)-\hat{v}_{\rho}(W_{i})\alpha(X_{i})^{2}\right\}
+\Lambda_{r}(\alpha),\label{alpha obj reg}%
\end{equation}
for $\mathcal{A}_{n}\subset\Gamma$, where $\hat{v}_\rho(W)$ is an estimator
\vsedit{of $v_\rho(X)$, in the sense that:
\begin{align}\label{eqn:projected-metric}
    \|\hat{v}_\rho - v_\rho\|_{X}^2 := \E\left[\left(\E[\hat{v}_\rho(W)\mid X] - v_\rho(X)\right)^2\right] = o_p(1).
\end{align}}
We refer to this $\hat\alpha_\ell$ as a weighted Riesz regression estimator. When $v_\rho(x)$ is known, we can use $\hat{v}_\rho(W)=v_\rho(X)$, and the expectation of this objective function is
\eqref{eq:wtrieszreg}, plus a penalty. Steps (3) and (4) are modified accordingly
to:
\begin{align*}
\hat{\theta} &  =\frac{1}{n}\sum_{\ell=1}^{L}\sum_{i\in I_{\ell}}\left\{
m(W_{i},\hat{\gamma}_{\ell})+\hat{\alpha}_{\ell}(X_{i})\rho(W_{i},\hat{\gamma
}_{\ell})\right\}  \\
\hat{V} &  =\frac{1}{n}\sum_{\ell=1}^{L}\sum_{i\in I_{\ell}}\left\{
m(W_{i},\hat{\gamma}_{\ell})-\hat{\theta}+\hat{\alpha}_{\ell}(X_{i})\rho
(W_{i},\hat{\gamma}_{\ell})\right\}  ^{2}.
\end{align*}

\begin{example}[Inverse Propensity Score Weighting\label{ex:ips}]\
The propensity score is useful for recovering counterfacutual distributions from observational data by weighting using the inverse propensity score (\cite{horvitz1952generalization}).
The superior performance of the automatic debiased machine learner in \cite{chernozhukov2022riesznet}, which is based on estimating the inverse of the propensity score directly, suggests the potential usefulness of this approach more generally. In this example we consider estimators of counterfactual averages based on estimators of the inverse propensity score.

To describe the estimators let $D$ be a treatment indicator, $Y$ be an outcome 
variable, with counterfactual value $Y(1)$ satisfying $Y(1)D=YD,$ $Z$ be
covariates, and $\gamma_{0}(Z)=1/\Pr(D=1\mid Z).$ be the inverse propensity score.
When $D$ and $Y(1)$ are independent conditional on $Z$ and $\gamma_{0}(Z)$ is
finite with probability one, the mean $\theta_{0}$ of $Y(1)$ is given by%
\[
\theta_{0}=\E[m(W,\gamma_{0})],\text{ }m(W,\gamma)=DY\gamma(Z).
\]
Also, the inverse propensity score satisfies%
\[ \E[\rho(W,\gamma_{0})\mid Z]=0,\text{ }\rho(W,\gamma)=1-D\gamma(Z).
\]
This conditional moment restriction can be interpreted as balancing for all possible functions of the covariates. This means that $\gamma_{0}(Z)$ is a generalized regression where $\Gamma$ is all functions of $Z$ with finite second moment
and the residual is $\rho(W,\gamma)=1-D\gamma(Z)$. Furthermore, the
conditional moment restriction corresponds to the first order condition for %
\[
\gamma_{0}=\arg\min_{\gamma} \E[-2\gamma(Z)+D\gamma(Z)^{2}].
\]
Thus $\gamma_{0}$ can be estimated by minimizing the sample average of
$-2\gamma(Z)+D\gamma(Z)^{2}$ over some set $\Gamma_n$ of functions of $Z$, as in
\[
\hat{\gamma}_{\ell}=\arg\min_{\gamma\in\Gamma_n}\left[  \sum_{i\notin I_{\ell}%
}\{-2\gamma(Z_{i})+D_{i}\gamma(Z_{i})^{2}\}+\Lambda_{r_{\gamma}}%
(\gamma)\right]  ,
\]
Also here $v_{\rho}(Z)=-\E[D\mid Z],$ so that we can take $v_{\rho}(W)=-D,$ and obtain $\hat{\alpha}_{\ell}$ as%
\begin{align*}
\hat{\alpha}_{\ell}  & =\arg\min_{\alpha\in\mathcal{A}_{n}}[\sum_{i\notin
I_{\ell}}\{-2D_{i}Y_{i}\alpha(Z_{i})+D_{i}\alpha(Z_{i})^{2}%
\}+\Lambda_{r}(\alpha)]\\
& =\arg\min_{\alpha\in\mathcal{A}_{n}}[\sum_{i\notin I_{\ell}}D_{i}\{
{Y}_{i}-\alpha(Z_{i})\}^{2}+\Lambda_{r}(\alpha)],
\end{align*}
where the last equality follows by adding $D_i{Y_i}^2$ inside the brackets, which does not affect the minimizer, and completing the square. Here we see that 
$\hat\alpha$ is a least squares learner of $\E[Y \mid D=1,Z]$. The resulting estimator of the parameter of interest is
\begin{align*}
\hat{\theta} &  =\frac{1}{n}\sum_{\ell=1}^{L}\sum_{i\in I_{\ell}}\left\{
D_iY_i\hat{\gamma}_{\ell}(Z_i)+\hat{\alpha}_{\ell}(Z_{i})[1-D_i\hat{\gamma
}_{\ell}(Z_i)]\right\}  \\
&
=\frac{1}{n}\sum_{\ell=1}^{L}\sum_{i\in I_{\ell}}\left\{
\hat{\alpha}_{\ell}(Z_{i})+D_i\hat{\gamma
}_{\ell}(Z_i)[Y_i-\hat{\alpha}_{\ell}(Z_{i})]\right\}.
\end{align*}
Here $\hat{\theta}$ has the classic doubly robust form \cite{robins1995semiparametric} of an average regression plus a bias correction term, with the key feature that the estimator $\hat{\gamma}_{\ell}(Z_i)$ of the inverse of the propensity score appears in place of the inverse of a propensity score estimator.
\end{example}

Below we give regularity conditions and a theorem to extend the results of subsection \ref{sec:largesampalpha} to the case where $v_\rho(x)$ is unknown and needs to be estimated. \vsedit{For simplicity of exposition we will only consider the case where the estimator is defined over a growing sieve space $\mathcal{A}_n$ and no regularization is used, i.e. $\Lambda_r(\alpha)=0$:
\begin{equation}
\hat{\alpha}=\arg\min_{\alpha\in\mathcal{A}_{n}}\sum_{i=1}^n \left\{  -2m(W_{i},\alpha)-\hat{v}_{\rho}(W_{i})\alpha(X_{i})^{2}\right\},\label{alpha obj}%
\end{equation}
Our estimation rate can easily be extended to regularized estimation with appropriate regularization weight.} Let $\mathcal{V}_n$ denote the function space in which the estimator $\hat{v}_\rho$ is restricted to lie in. Let $\alpha_{\ast}$ be any function in $\C{A}_n$ (e.g. we will typically consider $\alpha_{\ast}=\inf_{\alpha\in \C{A}_n} \|\alpha_{\ast}-\alpha_0\|$, but $\alpha_{\ast}$ can in fact be any function that is not chosen based on the samples). Define:
\vsedit{\begin{align*}
{\mathrm{star}}(\sqrt{\mathcal{V}_n} \cdot (\mathcal{A}_n-\alpha_{*}))=~  &  \{w\rightarrow
\xi\,\sqrt{|v(w)|}\,(\alpha(x)-\alpha_{*}(x)):\alpha%
\in\mathcal{A}_n, v\in \mathcal{V}_n,\text{ }\xi\in\lbrack0,1]\}\\
{\mathrm{star}}(m\circ\mathcal{A}_n-m\circ\alpha_{*})=~  &
\{w\rightarrow\xi\,(m(W,\alpha)-m(W,\alpha_{*}%
)):\alpha\in\mathcal{A}_n,\text{ }\xi\in\lbrack0,1]\}
\end{align*}}
\begin{assumption}\label{ass2bis}
$\left\Vert f\right\Vert _{\infty}\leq1$ for all $f\in{\mathrm{star}}(\sqrt{\mathcal{V}_n} \cdot (\mathcal{A}_n-\alpha_{*}))$ and $f\in {\mathrm{star}}(m\circ\mathcal{A}_n-m\circ\alpha_{*})$.
\end{assumption}
We remark that the uniform upper bound of $1$ can be replaced by any constant upper bound $b$ and the rate that we achieve will be identical, up to an extra multiplicative factor $b$, via a standard re-scaling argument (i.e. applying our result to re-scaled version of the original problem and then scaling back the guarantee).
\begin{assumption}\label{ass3bis}
The function $v_\rho$ and its estimate $\hat{v}_\rho$ satisfy that $\vsedit{|\hat{v}_\rho(W)|}, |v_\rho(X)| \leq C$, almost surely, and that for any $\alpha \in \C{A}_n$, the true function $v_\rho$ satisfies:
\begin{equation}
-\E[v_\rho(X) (\alpha(X)-\alpha_{\ast}(X))^2] \geq \lambda \E[(\alpha(X)-\alpha_{\ast}(X))^2], \label{eqn:strong-conv}
\end{equation}
for some constants $\lambda, C>0$. For notational convenience, $\lambda \leq 1$.
\end{assumption}

\begin{theorem}\label{thm1vunknown}
Let $\delta_{n}$ be an upper bound on the critical radius of $\mathrm{star}(\sqrt{\mathcal{V}_n}\cdot (\mathcal{A}_n-\alpha_{*}))$ and $\mathrm{star}(m\circ \mathcal{A}_n-m\circ\alpha_{*})$. If Assumptions \ref{ass1}, \ref{ass2bis} and \ref{ass3bis} are
satisfied then it follows that with probability $1-\zeta$, for some universal constant $C$,
\[
\left\Vert \hat{\alpha}-\alpha_{0}\right\Vert ^{2}\leq C\left(\frac{M}{\lambda^2}\delta_{n}%
^{2}+ \frac{1}{\lambda} \left\Vert \alpha_{\ast}-\alpha_{0}\right\Vert ^{2} + \frac{1}{\lambda^2} \vsedit{\|\hat{v}_\rho - v_\rho\|_X^2} + \frac{M\ln(1/\zeta
)}{\lambda n}\right).
\]
\end{theorem}
Moreover, we note that if a separate sample was used to estimate $\hat{v}_{\rho}$ and not the same as the one that was used for $\hat{\alpha}$, then we can weaken the theorem to only require $\delta_n$ to upper bound the critical radius of $\mathrm{star}(\mathcal{A}_n-\alpha_{*})$ and not $\mathrm{star}(\sqrt{\mathcal{V}_n}\cdot (\mathcal{A}_n-\alpha_{*}))$. Note that a sufficient condition for Equation~\eqref{eqn:strong-conv} is that $|v_\rho(X)|\geq \lambda$, almost surely. However, for most function spaces $\C{A}_n$, this condition can be satisfied by more benign assumptions. For this it is crucial that we only invoked the property at the difference of two functions that both lie in $\C{A}_n$ and not for instance for $\alpha - \alpha_0$ (since $\alpha_0$ can potentially lie outside of the space). For instance, if the functions in $\C{A}_n$ are are high-dimensional linear functions $\phi(X)'\beta$, then Equation~\eqref{eqn:strong-conv} is satisfied, with $\lambda = \mu/C$, if:
\begin{align*}
    \E[|v_\rho(X)|\phi(X)\phi(X)'] \succeq~& \mu I, &
    \E[\phi(X)\phi(X)'] \preceq~& C I,
\end{align*}
since then, if we let $\alpha=\phi(\cdot)'\beta$ and $\alpha_{\ast}=\phi(\cdot)'\beta_{\ast}$ and $\nu=\beta - \beta_{\ast}$, then:
\begin{align*}
    \|\alpha - \alpha_{\ast}\|^2 =~& \nu' \E[\phi(X)\phi(X)']\nu \leq C \|\nu\|_2^2\\
    \leq~& \frac{C}{\mu} \nu'\E[|v_\rho(X)| \phi(X)\phi(X)']\nu = \frac{C}{\mu}\E[|v_\rho(X)| (\alpha(X)-\alpha_{\ast}(X))^2] 
\end{align*}

The following assumption provides regularity conditions on the residual $\rho(w, \gamma)$ and the functional of interest $m(w, \gamma)$ for the generalized regression case.

\begin{assumption}
    \label{ass5} (i)
$\alpha_0(X)$ and $\E[\rho(W, \gamma_0)^2 \mid X]$ are bounded and $\E[m(W, \gamma_0)^2] < \infty$; (ii) $\hat{\alpha}_{\ell}(X)$ is bounded; (iii) $\E[\{\rho(W, \gamma) - \rho(W, \gamma_0)\}^2] \to 0$ if $\Vert \gamma - \gamma_0\Vert \to 0$; (iv) there is $C > 0$ such that for all $\Vert \gamma - \gamma_0\Vert$ small enough, $\E[\{\Bar{\rho}(W, \gamma) - \Bar{\rho}(W, \gamma_0)\}^2] \leq C \Vert \gamma - \gamma_0\Vert^2$, where $\Bar{\rho}(X, \gamma) = \E[\rho(W, \gamma) \mid X]$.
\end{assumption}

 The next condition allows for
$\rho(W,\gamma)$ to be nonlinear in $\gamma.$

\begin{assumption}
    \label{ass6}
    Either $\rho(W,\gamma)$ is affine in
$\gamma$ or $n^{1/4}\left\Vert \hat{\gamma}_{\ell}-\gamma
_{0}\right\Vert \overset{p}{\to}0$ and there are
$C, \varepsilon >0$ such that
\[
\left\vert \E[m(W,\gamma)-\theta_{0}+\alpha_{0}(X)\rho(W,\gamma)]\right\vert
\leq C\left\Vert \gamma-\gamma_{0}\right\Vert ^{2}.
\]
whenever $\left\Vert \gamma-\gamma_{0}\right\Vert ^{2} \leq \varepsilon$.
\end{assumption}

This assumption imposes the usual faster than $n^{-1/4}$ convergence rate for $\hat{\gamma}_{\ell}$ when $\rho(w,\gamma)$ is nonlinear in $\gamma$ but does not require that rate when $\rho(W,\gamma)$ is linear in $\gamma$. 

We have the following large sample inference result under these conditions.

\begin{theorem}\label{thm2vunknown}
    If Assumptions \ref{ass1}, \ref{ass4}, \ref{ass5} and \ref{ass6} are satisfied, then
    \[
    \sqrt{n}(\hat{\theta}-\theta_{0})\overset{p}{\to}N(0,V)\quad\text{and}\quad\hat{V}\overset{p}{\to}V.
    \]
where $V=\E[\{m(W,\gamma_{0})-\theta_{0}+\alpha_{0}(X)\rho(W, \gamma_0)\}^{2}]$.
\end{theorem}

\subsection{Nonlinear Effects of Multiple Regressions}\label{sec:multip}
Some important objects of interest are expectations of nonlinear functionals of multiple regressions. In this Section we give Auto-DML for such effects. Such effects have the form $\theta_{0}=\E[m(W,\gamma_{0})]$ where $m(w,\gamma)$ is nonlinear in a possible value $\gamma$ of multiple generalized regressions $(\gamma_{1}(X_{1}),...,\gamma_{J}(X_{J}))^{\prime}$ with regressors $X_{j},$ residual $\rho_{j}(W,\gamma_{j})$, and $\Gamma_{j}$ specific to each regression $\gamma_{j}(X_{j})$. The corresponding orthogonal score like is like that of subsection \ref{sec:gener} except that the bias correction is a sum of $J$ terms with the $j^{th}$ term being the bias correction for the learner of $\gamma_{j}$. Similarly to \cite{newey1994asymptotic}, pg. 1357, the orthogonal score is
\begin{equation}
\psi(w,\gamma,\alpha,\theta)=m(W,\gamma)-\theta+\sum_{j=1}^{J}\alpha_j(X_j)\rho_j(W,\gamma_j),\text{
}\gamma_j,\alpha_j\in\Gamma_j.\label{eq:orthmomnonlineargenreg}%
\end{equation}

\vscomment{We should explain here using the same "first-order" correction intuition, why we are using the directional derivative as the moment. Currently we don't offer any intuition}

Each of the terms in the bias correction can be estimated by the product of a learner $\hat{\alpha}_{j\ell}(X_{j})$ and the residual $\rho_{j}(W,\hat{\gamma}_{j\ell})$, but now the learner $\hat{\alpha}_{j\ell}(X_{j})$ differs from the one given in section \ref{sec:gener} in the way needed to account correctly for nonlinearity of $m(W,\gamma)$ in $\gamma$. The difference is that in the objective function for $\hat{\alpha}_{j\ell}(X_{j})$ the functional of interest $m(w,\alpha)$ is replaced by an estimated Gateaux derivative with respect to the $j^{th}$ component of $\gamma.$ Let
\[
\hat{D}_{j}(W,\alpha_{j})=\left.  \frac{d}{d\tau}m(W,\hat{\gamma
}_{\ell}+\tau e_{j}\alpha_{j})\right\vert _{\tau=0}%
\]
be such a Gateaux derivative estimator, where $e_{j}$ denotes the $j$-th column of the identity matrix. This derivative will often be straightforward to calculate as an analytic derivative with respect to the scalar $\tau.$ When $m(w,\gamma)$ is linear in a single $\gamma$ this derivative just evaluates $m(W_{i}, \gamma)$ at $\gamma=\alpha$ giving the $m(W,\alpha)$ of subsection \ref{sec:estim}.

To obtain $\hat{\alpha}_{j\ell}(X_{j})$ we also make use of an estimated derivative $\hat{v}_{\rho j}(W_{i})$ of $\rho_{j}(W,\gamma_{j})$ with respect to $\gamma_{j}$ at $\hat{\gamma}_{j\ell}.$ Then $\hat{\alpha}_{j\ell}$ is
given by
\begin{equation}
 \hat{\alpha}_{j\ell}= \arg\min_{\alpha_{j}\mathcal{\in A}^{j}_n}\left\{  \sum_{i\notin I_{\ell}}[-2\hat{D}_{j}%
(W_{i},\alpha_{j})-\hat{v}_{\rho j}(W_{i})\alpha_{j}(X_{ji})^{2}]\right\}
,   \label{eq:ahatnonl}
\end{equation}
where $\mathcal{A}^{j}_n$ is the set of approximating functions for $\alpha_{j}.$ As with linear $m(w, \gamma)$ this $\hat{\alpha}_{j\ell}$ depends just on $m(w,\gamma)$ and the first step. Thus $\hat{\alpha}_{j\ell}$ is automatic, in the same way as in section \ref{sec:linear}, in only requiring $m(w,\gamma)$ and the regression residual $\rho_{j}(W_{i},\gamma_{j})$ for its construction.

Below we give two examples of this setting:

\begin{example}[Marginal Effect in a Generalized Regression Model]
    Suppose that $X = (D,Z)$, where $D$ is a continuous treatment or policy variable and $Z$ are covariates. We are interested in $\theta_0 = \E[m(W, \gamma_0)]$ with $$m(W, \gamma) = \partial_a \mu (\gamma(X)) \partial_d \gamma(X).$$ The function $\gamma_0 \in \Gamma$ is assumed to satisfy the orthogonality condition \eqref{eq:orthotresid} for $\rho(W,\gamma)=\lambda(\gamma(X))[Y-\mu(\gamma(X))]$. This is the first order condition of a Generalized Regression Model with link function $\mu(a)$ (\cite{nelder1972generalized}). For example, when $Y$ is binary, $\mu(a)$ is a CDF and $\lambda(a) = \partial_a\mu(a)/[\mu(a)(1-\mu(a))]$, and $\Gamma$ is the mean square closure of finite linear combinations of $X$, this corresponds to a high dimensional, approximately sparse binary response model.

    In this example, $J = 1$, $\hat D(W, \alpha) = \partial_a^2(\hat\gamma(X))\partial_d \hat\gamma(X) \alpha(X) + \partial_a \mu(\hat\gamma(X))\partial_d \alpha (X)$ and $\hat v_\rho(X) = \partial_a\lambda(\hat\gamma(X))[Y - \mu(\hat\gamma(X))] - \lambda(\hat\gamma(X))\partial_a\mu(\hat\gamma(X))$. In the Logit case, $\lambda(a) \equiv 1$, and so this simplifies to $\hat v_\rho(X) = -\partial_a\mu(\hat\gamma(X))$. The weighted Riesz regression estimator  $\hat\alpha_\ell$ can be found by evaluating \eqref{eq:ahatnonl} at these quantities, and then used to build the Neyman orthogonal score \eqref{eq:orthmomnonlineargenreg}.
\end{example}

\begin{example}[Inverse Logit Propensity Score Weighting]
    Suppose now that $Y$ is a continuous or discrete outcome, $X = (D, Z)$  where $D$ is a binary treatment and $Z$ are covariates, and the parameter of interest is $\theta_0 = \E[m(W, \gamma_0)]$ with $$m(W, \gamma) = \frac{DY}{\Lambda(\gamma(Z))},$$ where $\Lambda(a)$ is the standard logistic CDF. The parameter $\gamma_0 \in \Gamma$ satisfies the orthogonality condition \eqref{eq:orthotresid} for $\rho(W,\gamma)=D-\Lambda(\gamma(Z))$. This corresponds to inverse propensity score weighting of the outcome, where the propensity score is modelled by a flexible Logit specification. For example, if $\Gamma$ is the mean square closure of finite linear combinations of $X$, this corresponds to a high dimensional, approximately sparse logit model; if $\Gamma$ is the space of all square-integrable functions, the model is essentially unrestricted.

    In this example, $J = 1$, $\hat D(W, \alpha) = - DY \partial_a\Lambda(\hat\gamma(X)) \alpha(X)/\Lambda(\hat\gamma(X))^2$ and $\hat v_\rho(X) = -\partial_a\Lambda(\hat\gamma(X))$. The debiasing function $\hat\alpha_\ell$ can be found by evaluating \eqref{eq:ahatnonl} at these quantities, and then used to build the Neyman orthogonal score \eqref{eq:orthmomnonlineargenreg}.
\end{example}

It is straightforward to obtain a convergence rate for $\hat{\alpha}_{j\ell}$ analogous to Theorem \ref{thm1vunknown}. The following result does so while accounting for the presence of $\hat{\gamma}$ in $\hat{D}_{j}(W_{i},\alpha_{j}).$ For notational
convenience we suppress the $j$ subscripts.

\begin{assumption}\label{ass-der-estimate}
The estimate $\hat{D}$ satisfies that:
\begin{align}
    |\E[\hat{D}(W, \alpha) - D(W, \alpha)]| \leq \epsilon_{mn}\, \|\alpha\|
\end{align}
\end{assumption}
\begin{theorem} \label{thm1nonlin} If the conditions of Theorem \ref{thm1vunknown} and Assumption~\ref{ass-der-estimate} is satisfied then it follows that with probability $1-\zeta$, for some universal constant $C$,
\[
\left\Vert \hat{\alpha}-\alpha_{0}\right\Vert ^{2}\leq C\left(\frac{M}{\lambda^2}\delta_{n}%
^{2}+ \frac{1}{\lambda} \left\Vert \alpha_{\ast}-\alpha_{0}\right\Vert ^{2} + \frac{1}{\lambda^2} \left(\vsedit{\|\hat{v}_\rho - v_\rho\|_X^2} + \epsilon_{mn}^2\right) + \frac{M\ln(1/\zeta
)}{\lambda n}\right).
\]
\end{theorem}

The construction of $\hat{\theta}$ is analogous to that in subsection \ref{sec:estim} with the bias correction term being the sum of terms for each $\gamma_{j}$ in $\gamma.$ That is,
\begin{align}
\hat{\theta}  &  = \frac{1}{n}\sum_{\ell=1}^{L}\sum_{i\in I_{\ell}}%
\{m(W_{i},\hat{\gamma}_{\ell})+\sum_{j=1}^{J}\hat{\alpha}_{j\ell}(X_{ji}%
)\rho_{j}(W_{i},\hat{\gamma}_{j\ell})\},\label{nonlin est}\\
\hat{V}  &  =\frac{1}{n}\sum_{\ell=1}^{L}\sum_{i\in I_{\ell}} \left\lbrace m(W_{i},\hat{\gamma}_{\ell}%
)-\hat{\theta}+\sum_{j=1}^{J}\hat{\alpha}_{j\ell}(X_{ji})\rho_{j}(W_{i}%
,\hat{\gamma}_{j\ell}) \right\rbrace^2.\nonumber
\end{align}

It is straightforward to specify conditions for asymptotic normality of $\hat{\theta}$ and consistency of $\hat{V}$ by combining the conditions of section \ref{sec:gener} with the convergence rate result of Corollary \ref{cor1}. For relative simplicity we give a result only for neural net learners. We also assume for simplicity that each $X_{j}$ has the same dimension $d.$

\begin{assumption}\label{ass7}
$\E[m(W,\gamma_{0})^{2}]<\infty$ and for each
$j,$ (i) $\E[\rho_{j}(W,\gamma_{j0})^{2}\mid X]$ is bounded (ii),
$\E[\{\rho_{j}(W,\gamma_{j})-\rho_{j}(W,\gamma_{0})\}^{2}]\to0$
if $\left\Vert \gamma_{j}-\gamma_{j0}\right\Vert \to0$;
(iii) there is $C>0$ such that for all $\left\Vert \gamma
_{j}-\gamma_{j0}\right\Vert $ small enough $\E[\{\bar{\rho}_{j}%
(X_{j},\gamma_{j})-\bar{\rho}_{j}(X,\gamma_{j0})\}^{2}]\leq C\left\Vert\gamma_{j}-\gamma_{j0}\right\Vert ^{2}$, where $\bar{\rho}_{j}(X_{j},\gamma_{j})=\E[\rho_{j}(W,\gamma_{j})\mid X_{j}].$    
\end{assumption}
This condition is analogous to Assumption \ref{ass5}. 

\begin{assumption}\label{ass8}
$n^{1/4}\left\Vert \hat{\gamma}_{j\ell}-\gamma_{j0}\right\Vert \overset{p}{\to}0$  for each $j$ and there are $C,\varepsilon>0$ such that for
\[
\left\vert \E[m(W,\gamma)-\theta_{0}+\sum_{j=1}^{J}\alpha_{j0}(X_{j})\rho(W,\gamma_{j})]\right\vert \leq C\sum_{j=1}^{J}\left\Vert \gamma_{j}-\gamma_{j0}\right\Vert ^{2}.
\]
whenever $\left\Vert \gamma_{j}-\gamma_{j0}\right\Vert ^{2} \leq \varepsilon$ for all $j = 1, \ldots, J$.
\end{assumption}
This condition is analogous to Assumption \ref{ass6}.

\begin{theorem} \label{thm2nonlin}
    If Assumptions \ref{ass1}, \ref{ass4}, \ref{ass7} and \ref{ass8} are satisfied 
    for each $j = 1, \ldots, J$, then
\[
\sqrt{n}(\hat{\theta}-\theta_{0})\overset{p}{\to}N(0,V),\text{
}\hat{V}\overset{p}{\to}V.
\]
where $V=\E[\{m(W,\gamma_{0})-\theta_{0}+\sum_{j=1}^{J}\alpha_{j0}(X_{j}%
)\rho_{j}(W,\gamma_{j0})\}^{2}]$.
\end{theorem}

\section{Empirical Application \label{sec:empir}}
To illustrate our methods, we study whether applicant race is a significant predictor of banks' mortgage denial decisions. Following \cite{munnell1996mortgage}, we use the publicly available Boston Home Mortgage Disclosure Act (HDMA) dataset. The dataset contains information on 2,925 mortgage applications made in 1990 in the Greater Boston metropolitan area. We restrict attention to black and white applicants in single-family households (excluding other racial minorities and multi-family residences), which reduces our sample size to 2,380 observations.

Our outcome of interest is an indicator $Y = 1$ if the mortgage application was denied. Our regressor of interest is an indicator $D = 1$ if the applicant is black. We also have access to a vector of covariates, which we denote by $Z$, containing financial and other characteristics of the applicant that banks may factor into their mortgage denial decisions. These include monthly debt to income (DTI) ratio; monthly housing expenses to income (HTI) ratio; loan to assessed property value (LTV) ratio; a categorical variable for ``bad'' consumer credit score with 6 categories (1 if no slow payments or delinquencies, 2 if one or two slow payments or delinquencies, 3 if more than two slow payments or delinquencies, 4 if insufficient credit history for determination, 5 if delinquent credit history with payments 60 days overdue, and 6 if delinquent credit history with payments 90 days overdue); a categorical variable for ``bad'' mortgage credit score with 4 categories (1 if no late mortgage payments, 2 if no mortgage payment history, 3 if one or two late mortgage payments, and 4 if more than two late mortgage payments); an indicator for public record of credit problems including bankruptcy, charge-offs, and collective actions; an indicator for denial of application for mortgage insurance; three indicators for self-employed, single, and high school graduate, the 1989 Massachusetts unemployment rate in the applicant's industry, and an indicator for whether the unit is a condominium.

\begin{table}[] \centering
    \caption{Summary Statistics}
    \label{tab:summary}
    \begin{tabular}{rcccccc}
    \toprule
    & \multicolumn{2}{c}{Full Sample} & \multicolumn{2}{c}{Black} & \multicolumn{2}{c}{White} \\
    \cmidrule(r){2-3} \cmidrule(lr){4-5} \cmidrule(l){6-7}
    & mean & sd & mean & sd & mean & sd \\
    \midrule
    Deny & 0.12 & 0.32 & 0.28 & 0.45 & 0.09 & 0.29 \\
Monthly DTI Ratio & 0.19 & 0.01 & 0.19 & 0.01 & 0.19 & 0.01 \\
Monthly HTI Ratio & 0.12 & 0.01 & 0.12 & 0.01 & 0.12 & 0.01 \\
LTV Ratio & 0.37 & 0.03 & 0.38 & 0.01 & 0.37 & 0.03 \\
Consumer Credit Ind. & 2.12 & 1.67 & 3.02 & 2.01 & 1.97 & 1.55 \\
Mortgage Credit Ind. & 1.72 & 0.54 & 1.88 & 0.42 & 1.69 & 0.55 \\
Public Record & 0.07 & 0.26 & 0.18 & 0.38 & 0.06 & 0.23 \\
Denied Insurance & 0.02 & 0.14 & 0.05 & 0.22 & 0.02 & 0.12 \\
Self-Employed & 0.12 & 0.32 & 0.07 & 0.26 & 0.12 & 0.33 \\
Single & 0.39 & 0.49 & 0.52 & 0.50 & 0.37 & 0.48 \\
High School & 0.98 & 0.13 & 0.97 & 0.18 & 0.99 & 0.12 \\
Industry Unemp. & 3.77 & 2.03 & 3.45 & 1.50 & 3.83 & 2.10 \\
Condominium & 0.29 & 0.45 & 0.49 & 0.50 & 0.25 & 0.44 \\
\midrule
$N$ & \multicolumn{2}{c}{2,380} & \multicolumn{2}{c}{339} & \multicolumn{2}{c}{2,041} \\
\bottomrule
 
    \end{tabular}
\end{table}

Table \ref{tab:summary} reports the sample means and standard deviations of the variables used in the analysis. The probability of being denied a mortgage is 19 percentage points higher for black applicants than for white applicants. However, black applicants are also more likely to have financial and socio-economic characteristics linked to mortgage denial, as Table \ref{tab:summary} shows. For example, black applicants have higher (worse) consumer and mortgage credit indices on average, and are more likely to have a public record of credit problems and to be single. We would like to test whether the racial differences in probability of mortgage denial persist once we control for these covariates.

To showcase the versatility of our method, we present results for three estimands:
\begin{enumerate}
    \item \textit{Difference in Probability of Mortgage Denial:}
    \[\theta_0 = \E[\gamma_0(1, Z) - \gamma_0(0, Z)], \quad \text{where} \quad \gamma_0(D, Z) = \Pr(Y = 1 \mid D, Z) = \E[Y \mid D, Z].\]
    This is an average linear effect for a conditional mean (Section \ref{sec:linear}). This parameter can be interpreted as an average difference in probability of mortgage denial between a black and a white applicant with the same value of covariates $Z$.

    \vspace{0.5em}
    \item \textit{Average Difference in Log-Odds of Mortgage Denial:}
    \[\theta_0 = \E[\gamma_0(1, Z) - \gamma_0(0, Z)], \quad \text{where} \quad \gamma_0(D, Z) = \ln \frac{\Pr(Y = 1 \mid D, Z)}{\Pr(Y = 0 \mid D, Z)}.\]
    This is an average non-linear effect for a generalized regression (Section \ref{sec:multip}). Because $\ln(a) - \ln(b) \approx (a-b)/b$ when $(a-b)/b$ is small, this parameter can be interpreted as an approximate average percentage difference in odds of mortgage denial between a black and a white applicant with the same value of covariates $Z$. As discussed in Section \ref{sec:multip}, this $\gamma_0$ minimizes the logistic regression loss function,
    \[\ell(W, \gamma) = Y \ln \Lambda(\gamma(X)) + (1 - Y) \ln [1 - \Lambda(\gamma(X))],\]
    with corresponding generalized residual
    \[\rho(W, \gamma) = Y - \Lambda(\gamma(X)),\]
    for $\Lambda(t) := 1/(1 + e^{-t})$, the standard logistic CDF.

    \vspace{0.5em}
    \item \textit{Average Difference in Odds of Mortgage Denial:}
    \[\theta_0 = \E[e^{\gamma_0(1, Z)} - e^{\gamma_0(0, Z)}], \quad \text{where} \quad \gamma_0(D, Z) = \ln \frac{\Pr(Y = 1 \mid D, Z)}{\Pr(Y = 0 \mid D, Z)}.\]
    This is an average non-linear effect for a generalized regression (Section \ref{sec:multip}). It can be interpreted as an average difference in odds of mortgage denial between a black and a white applicant with the same value of covariates $Z$.
\end{enumerate}

We estimate these parameters using AutoDML, where both $\hat{\gamma}$ and $\hat{\alpha}$ are neural net learners. For the difference in probability, we have $v_\rho(X) = -1$, since $\gamma_0$ is a conditional mean. For the average difference in log-odds and the average difference in odds, $v_\rho(X) = -\lambda(\gamma_0(X))$, where $\lambda(t) := e^{-t}/(1 - e^{-t})^2$ is the standard logistic PDF, which we estimate by replacing $\gamma_0$ with a preliminary estimate $\hat\gamma_{\text{prel}}$, also based on a neural net learner. We describe the architecture and training hyperparameter choice in detail in Appendix \ref{sec:hyperpar}.

\begin{table}[!h] \centering
    \caption{Empirical Application Results: Racial Differences in Probability, Average Log-Odds and Average Odds of Mortgage Denial}
    \label{tab:empirical}
    \begin{tabular}{rcccccc}
    \toprule
    & \multicolumn{2}{c}{Probability} & \multicolumn{2}{c}{Log-Odds} & \multicolumn{2}{c}{Odds} \\
    \cmidrule(r){2-3} \cmidrule(lr){4-5} \cmidrule(l){6-7}
    & est & se & est & se & est & se \\
    \midrule
    Main Spec. & 0.080 & (0.021) & 0.829 & (0.152) & 0.157 & (0.044) \\
 
    \bottomrule
    \end{tabular}
\end{table}

Table \ref{tab:empirical} presents the results of our main analysis. Once we control for covariates, the difference in probability of mortgage denial decreases from 19 to 8 percentage points. If we look at the average log-odds or odds instead, we observe differences of 0.829 or 0.157, respectively. These differences are all estimated to be statistically different from 0 at the 1\% significance level.

\begin{figure}
    \centering
    \subfloat[Differences in Probability by Consumer Credit Ind.]{
\includegraphics[width=.45\textwidth]{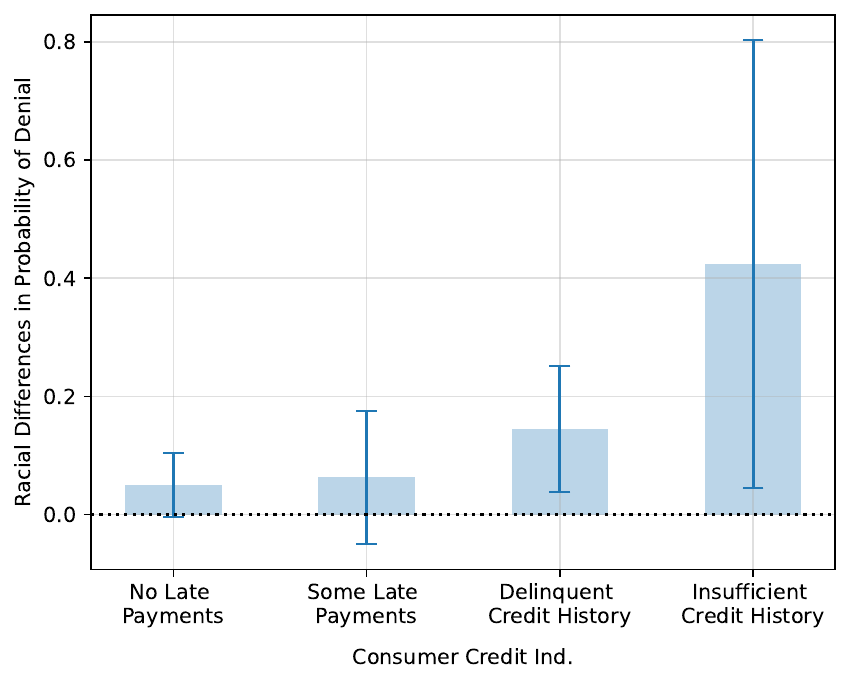}
    } \qquad
    \subfloat[Differences in Probability by Monthly DTI Ratio]{\includegraphics[width=.45\textwidth]{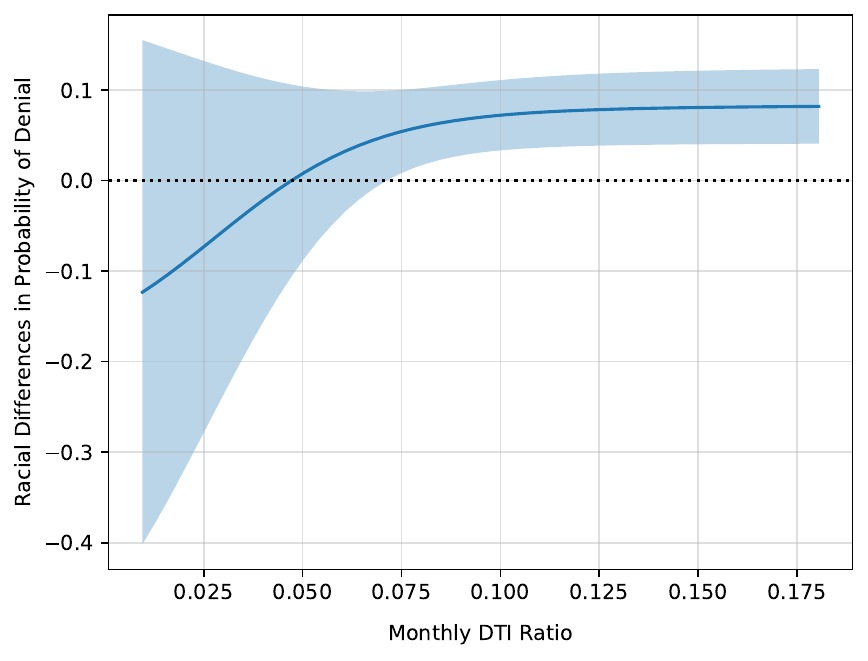}
    }

    \subfloat[Differences in Avg. Log-Odds by Consumer Credit Ind.]{
\includegraphics[width=.45\textwidth]{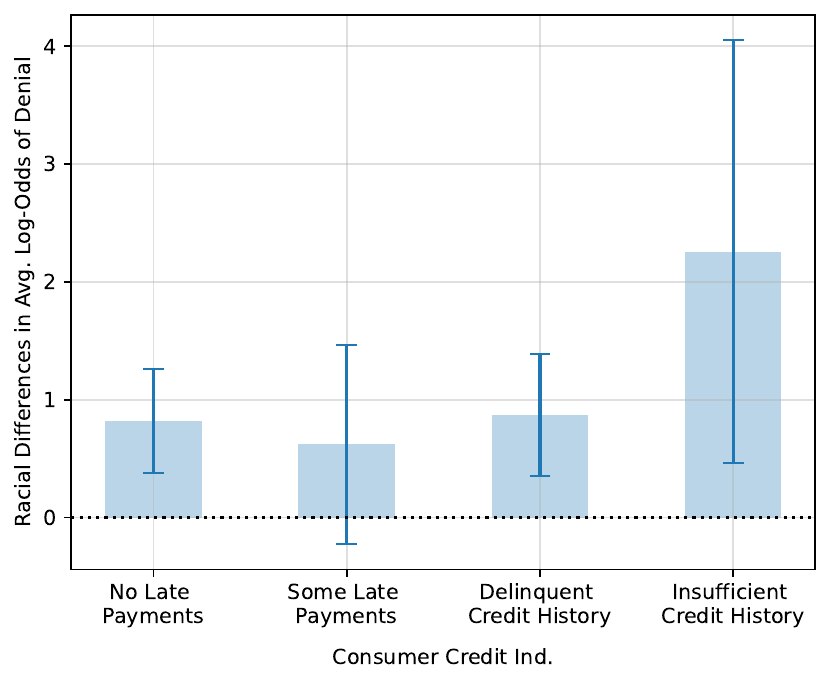}
    } \qquad
    \subfloat[Differences in Avg. Log-Odds by Monthly DTI Ratio]{\includegraphics[width=.45\textwidth]{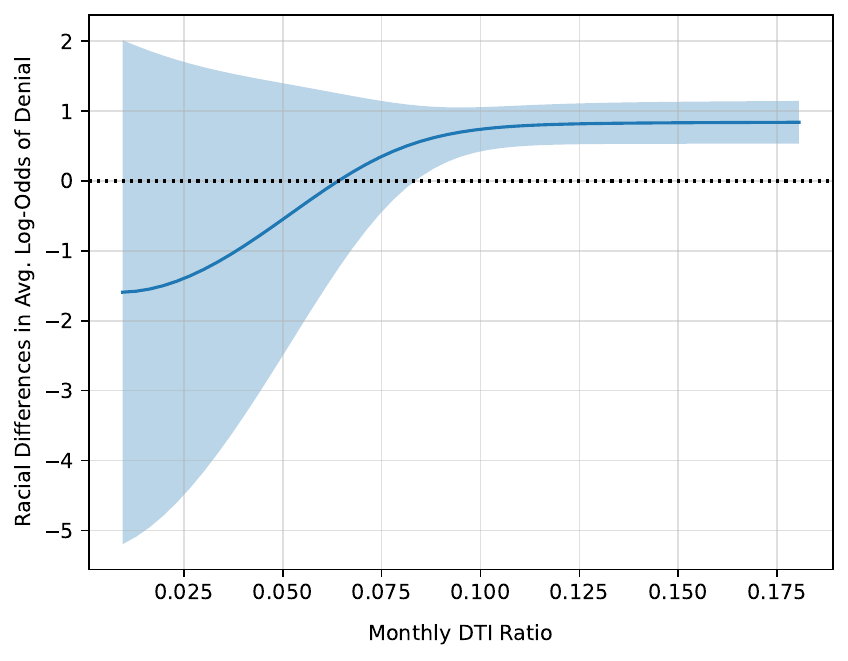}
    }

    \subfloat[Differences in Avg. Odds by Consumer Credit Ind.]{
\includegraphics[width=.45\textwidth]{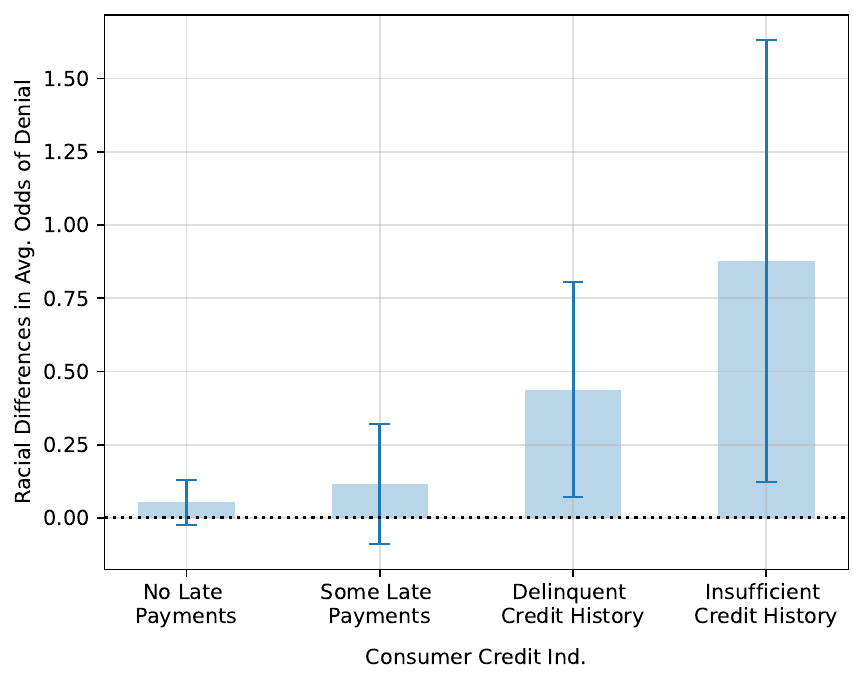}
    } \qquad
    \subfloat[Differences in Avg. Odds by Monthly DTI Ratio]{\includegraphics[width=.45\textwidth]{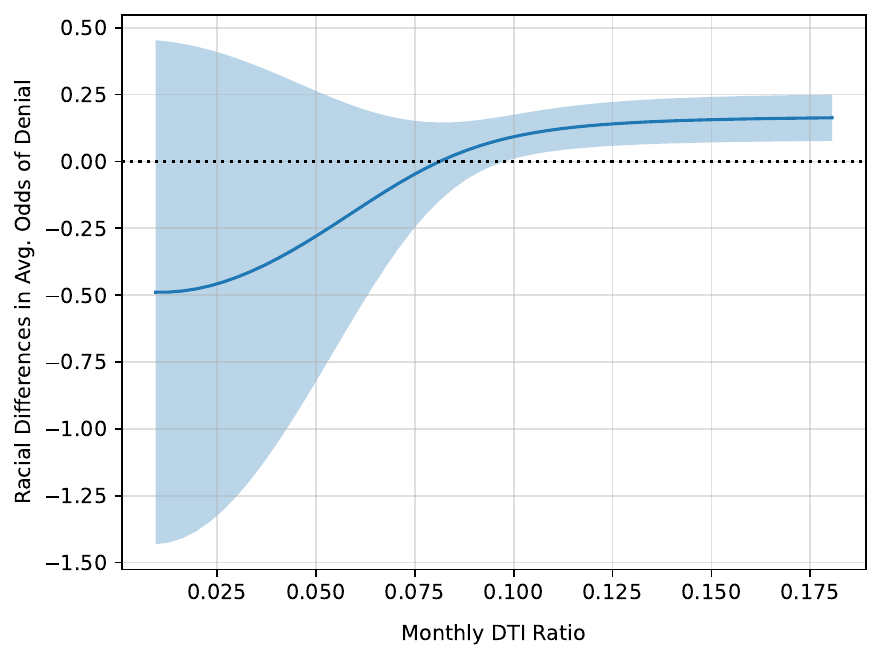}
    }
    
    \caption{Heterogeneous Effects}
    \label{fig:CATEs}
\end{figure}

A slight modification of our method allows us to estimate average differences for subgroups of applicants with certain characteristics (analogous to conditional average treatment effects or CATEs). Suppose we want to estimate an average effect for applicants with $Z_j = z$ for a particular covariate $Z_j$. To obtain these, we weight the Neyman orthogonal estimating equation \eqref{nonlin est} as follows:
\[\hat{\theta}(z) = \frac{1}{\sum_{i}^n \omega_i(z)}\sum_{\ell=1}^{L}\sum_{i\in I_{\ell}}%
\omega_i(z) \{m(W_
{i},\hat{\gamma}_{\ell})+\hat{\alpha}_{\ell}(X_{i}%
)\rho(W_{i},\hat{\gamma}_{\ell})\}.\]
When $Z_j$ is a categorical variable, we take $\omega_i(z) = 1\{Z_j = z\}$. When $Z_j$ is continuous, we take $\omega_i(z) = K((Z_i - z)/h)$ for a kernel function $K$ and a small but fixed bandwith $h$.\footnote{\cite{chernozhukov2019global} analyze a localized version of this parameter, that is, the limit as $h \to 0$, which is beyond the scope of this paper.} Figure \ref{fig:CATEs} presents the racial differences in probability, average log-odds and average odds by values of the consumer credit index and the monthly DTI ratio. Remarkably, we estimate the racial differences in all three estimands to be higher for applicants with a delinquent credit history or with insufficient credit history (although the latter is quite imprecisely estimated). The racial differences appear to be constant for most of the range of the monthly DTI ratio variable, except values below 0.075 for which it is also imprecisely estimated.

\section{Monte Carlo Simulations}\label{sec:simul}
\subsection{Based on the Empirical Application}\label{sec:simul1} First, we analyze the performance of our method in the setting of our empirical application. We redraw the covariates $Z$ based on a generative adversarial network (GAN) trained on the real mortgage data. We use an elastic-net Logit, with penalties chosen by cross-validation, to estimate the outcome regression $\Pr(Y = 1 \mid D, Z)$ and the propensity score $\Pr(D = 1 \mid Z)$ in the real mortgage data, which we take as ground truth in our simulations. We present results for the three estimands on interest in Section \ref{sec:empir}: the difference in probability, the average difference in log-odds and the average difference in odds.

Table \ref{tab:simul1} presents simulation results over 1,000 draws for $n = 2,000$ and $n = 10,000$. The first column shows the non-parametric $R^2$ for $\gamma$, defined as $R^2(\gamma) = 1 - \E[(\hat{\gamma}(X) - \gamma_0(X))^2]/\mathrm{Var}(\gamma_0(X))$, where the expectation is evaluated over a test set not used to estimate $\hat{\gamma}$. The second column shows the same non-parametric $R^2$ metric for $\alpha$. We also give the mean absolute error (MAE), bias, standard deviation (sd), the average standard error to standard deviation ratio (se/sd) and the coverage of a 95\% confidence interval (covg.).

For the estimands we study, the Riesz regression $\alpha_0$ can be characterized explicitly based on the Riesz representer for the Average Treatment Effect (ATE),
\[v_m(X) = \frac{D}{\Pr(D = 1 \mid Z)} - \frac{1-D}{1 - \Pr(D = 1 \mid Z)}.\]
For the difference in probability estimand, which is an average linear effect for a conditional mean (Section \ref{sec:linear}), $\alpha_0(X) = v_m(X)$. For the average difference in log-odds, which is an average linear effect for a generalized regression (Section \ref{sec:gener}), the weighted Riesz regression is $\alpha_0(X) = v_m(X)/(-v_\rho(X))$, where $v_\rho(X) = -\lambda(\gamma_0(X))$. Finally, for the average difference in odds, which is an average non-linear effect for a generalized regression (Section \ref{sec:multip}), we have $\alpha_0(X) = e^{\gamma_0(X)}v_m(X)/(-v_\rho(X))$. 

Our automatic debasing method does not make use of this explicit characterization of $\alpha_0$. To benchmark our results, we compare its performance to an estimator that uses the explicit characterization of $\alpha_0$. In the ATE setting, this is known as an Augmented Inverse Propensity Weighting (AIPW) estimator, so we will refer to these as AIPW-like. To build the AIPW-like estimator of $\theta_0$ we plug learners of the outcome and treatment propensities $\Pr(Y = 1 \mid D, Z)$ and $\Pr(D = 1 \mid Z)$ into the formula for $\alpha_0$; we try both a non-parametric version based on neural nets (the same architecture and hyperparameters as our main specification) and a ``well-specified'' version, where we use the same elastic-net Logit that we used to build the ground truth.

A comparison between auto-DML and the AIPW-like benchmark sheds light on the advantages of our automatic approach. For large sample sizes, $n = 10,000$, both methods perform comparably well. There is no loss in efficiency between our main specification (which uses a non-parametric, NN-based method) and the correctly-specified AIPW estimator (which uses correctly specified parametric learners for $\gamma$ and $\alpha$). Our automatic debiasing method achieves close to nominal coverage, whereas the AIPW-like method gets worse coverage when we use the non-parametric, neural net specification. A reason for that could be that, in the explicit characterization of $\alpha_0$, we are plugging numbers that are close to zero into a denominator (such as the propensity score $\Pr(D = 1 \mid Z)$ or the logit pdf $v_\rho(X)$), so that estimation error amplifies. This is reflected into the large negative non-parametric $R^2$ for $\alpha$ when $n = 2,000$. 

Consistent with these results, in work that followed up on the first version of this paper (\cite{chernozhukov2022riesznet}), we found that our automatic debiasing method using neural net and random forest Riesz regressions performed much better than state of the art methods based on inverse propensity score weighting in Monte Carlo experiments. \cite{singh2023double} also found that automatic debiased estimators performed better than plugin-based estimators in the setting of local average treatment effects.

\begin{table}[] \centering
    \caption{Simulation Results: Based on the Empirical Application}
    \label{tab:simul1}
    \begin{tabular}{r*{7}{c}}
    \toprule
    & \multicolumn{7}{c}{Probability} \\
    \cmidrule(r){2-8}
    & $R^2(\gamma)$ & $R^2(\alpha)$ & MAE & bias & sd & se/sd & covg. \\
    \midrule
    \multicolumn{8}{l}{$n = 2,000$} \\
    Main Spec. & 0.639 & 0.820 & 0.021 & 0.005 & 0.027 & 0.938 & 0.931 \\
 
    AIPW, NN & 0.622 & -1e7 & 0.022 & 0.005 & 0.028 & 0.890 & 0.923 \\
 
    AIPW, well spec. & 0.807 & 0.912 & 0.021 & 0.003 & 0.026 & 0.956 & 0.932 \\
 
    \addlinespace
    \multicolumn{8}{l}{$n = 10,000$} \\
    Main Spec. & 0.898 & 0.969 & 0.009 & 0.002 & 0.012 & 0.973 & 0.946 \\
 
    AIPW, NN & 0.892 & 0.969 & 0.009 & 0.002 & 0.012 & 0.970 & 0.938 \\
 
    AIPW, well spec. & 0.963 & 0.988 & 0.009 & 0.001 & 0.012 & 0.984 & 0.949 \\

    \addlinespace
    \midrule
    & \multicolumn{7}{c}{Log-Odds} \\
    \cmidrule(r){2-8}
    & $R^2(\gamma)$ & $R^2(\alpha)$ & MAE & bias & sd & se/sd & covg. \\
    \midrule
    \multicolumn{8}{l}{$n = 2,000$} \\
    Main Spec. & 0.583 & 0.368 & 0.203 & -0.010 & 0.265 & 0.871 & 0.899 \\
 
    AIPW, NN & 0.585 & -2e7 & 0.197 & 0.019 & 0.260 & 0.840 & 0.893 \\
 
    AIPW, well spec. &0.767 & -3e3 & 0.191 & -0.019 & 0.243 & 0.909 & 0.915 \\
 
    \addlinespace
    \multicolumn{8}{l}{$n = 10,000$} \\
    Main Spec. & 0.889 & 0.914 & 0.086 & 0.000 & 0.107 & 0.910 & 0.925 \\
 
    AIPW, NN & 0.888 & 0.925 & 0.084 & 0.007 & 0.104 & 0.927 & 0.928 \\
 
    AIPW, well spec. &0.959 & 0.959 & 0.083 & -0.007 & 0.104 & 0.977 & 0.938 \\

    \addlinespace
    \midrule
    & \multicolumn{7}{c}{Odds} \\
    \cmidrule(r){2-8}
    & $R^2(\gamma)$ & $R^2(\alpha)$ & MAE & bias & sd & se/sd & covg. \\ 
    \midrule
    \multicolumn{8}{l}{$n = 2,000$} \\
    Main Spec. & 0.589 & 0.653 & 0.039 & 0.015 & 0.055 & 0.786 & 0.932 \\
 
    AIPW, NN & 0.585 & -9e6 & 0.044 & 0.018 & 0.061 & 0.820 & 0.931 \\
 
    AIPW, well spec. &0.767 & -2e10 & 0.199 & 0.174 & 4.567 & 0.022 & 0.899 \\
 
    \addlinespace
    \multicolumn{8}{l}{$n = 10,000$} \\
    Main Spec. & 0.889 & 0.830 & 0.015 & 0.003 & 0.019 & 0.948 & 0.947 \\
 
    AIPW, NN & 0.888 & 0.867 & 0.016 & 0.006 & 0.020 & 0.960 & 0.936 \\
 
    AIPW, well spec. &0.959 & 0.866 & 0.017 & 0.005 & 0.021 & 0.994 & 0.95 \\
 
    \bottomrule
    \end{tabular}
\end{table}

\subsection{Additional Simulations} We present an additional set of simulations when $\gamma_0$ is a quantile of the conditional distribution of $Y \mid X$ for a continuous outcome $Y$. For simplicity, we will focus on the conditional median. As discussed in Section \ref{sec:gener}, this corresponds to the generalized residual $\rho(W, \gamma) = 0.5 - 1(Y < \gamma(X)),$ which is a sub-derivative of the ``check'' loss function $\ell(W, \gamma) = 0.5|Y - \gamma(X)|$.

We will consider four simulation settings. In the first setting, the object of interest will be the average difference in conditional median for a binary treatment $D$, that is, $\theta_0 = \E[\gamma_0(1, Z) - \gamma_0(0, Z)]$. We assume that $\Pr(D = 1 \mid Z) = \mathrm{logit}(-0.1 + 0.5 Z_1 - 0.2 Z_2)$, and that $Y = \mu_Y(X) + \epsilon$ with $\mu_Y(X) = 0.5D - 0.2 D \times Z_3 + 0.3 Z_2$ and $\epsilon \sim N(0, 1)$.  The second setting considers the same data generating process, but we use a twice-differentiable quantile loss based on Epanechnikov kernel smoothing, due to \cite{he2023smoothed}, rather than the non-differentiable check loss function. In the third setting, our object of interest is the average derivative of the conditional median with respect to a continuous treatment $D$, that is, $\theta_0 = \E[\partial_d \gamma_0(D, Z)]$, where we draw $D = \mu_D(Z) + \eta$ for $\mu_D(Z) = -0.1 + 0.5 Z_1 - 0.2 Z_2$ and $\eta \sim N(0, 1)$; the distribution of $Y \mid X$ is as before. The last setting uses the same data generating process, but it focuses on a non-linear parameter, the average derivative squared $\theta_0 = \E[(\partial_d \gamma_0(D, Z))^2]$. This parameter is useful, for example, in testing whether $\partial_d \gamma_0(D, Z) = 0$ with probability 1. In all settings, to work with a realistic data generating process for the covariates, we will draw $Z$ from the same GAN trained on the real mortgage data that we used in the previous subsection.

The results over 1,000 simulation draws are presented in Table \ref{tab:simul2}. In all the settings we consider, our estimator performs well, with low bias and coverage confidence intervals close to the nominal 95\% level.

\begin{table}[] \centering
    \caption{Simulation Results: Additional Designs}
    \label{tab:simul2}
    \begin{tabular}{r*{7}{c}}
    \toprule
    & \multicolumn{7}{c}{Quantiles} \\
    \cmidrule(r){2-8}
    & $R^2(\gamma)$ & $R^2(\alpha)$ & MAE & bias & sd & se/sd & covg. \\
    \midrule
    \multicolumn{8}{l}{$n = 2,000$} \\
    Main Spec. & 0.737 & 0.981 & 0.045 & -0.003 & 0.055 & 1.035 & 0.96 \\

    \addlinespace
    \multicolumn{8}{l}{$n = 10,000$} \\
    Main Spec. & 0.890 & 0.993 & 0.020 & 0.001 & 0.025 & 1.022 & 0.957 \\

    \addlinespace
    \midrule
    & \multicolumn{7}{c}{Quantiles, Smooth Loss} \\
    \cmidrule(r){2-8}
    & $R^2(\gamma)$ & $R^2(\alpha)$ & MAE & bias & sd & se/sd & covg. \\
    \midrule
    \multicolumn{8}{l}{$n = 2,000$} \\
    Main Spec. & 0.757 & -8.562 & 0.044 & -0.006 & 0.116 & 0.469 & 0.957 \\

    \addlinespace
    \multicolumn{8}{l}{$n = 10,000$} \\
    Main Spec. & 0.898 & 0.993 & 0.019 & 0.000 & 0.024 & 0.984 & 0.948 \\

    \addlinespace
    \midrule
        & \multicolumn{7}{c}{Quantiles, Smooth Loss, Continuous Treatment} \\
    \cmidrule(r){2-8}
    & $R^2(\gamma)$ & $R^2(\alpha)$ & MAE & bias & sd & se/sd & covg. \\
    \midrule
    \multicolumn{8}{l}{$n = 2,000$} \\
    Main Spec. & 0.869 & 0.944 & 0.021 & -0.002 & 0.026 & 0.968 & 0.928 \\

    \addlinespace
    \multicolumn{8}{l}{$n = 10,000$} \\
    Main Spec. & 0.952 & 0.979 & 0.010 & -0.002 & 0.012 & 0.941 & 0.922 \\

    \addlinespace
    \midrule
        & \multicolumn{7}{c}{Quantiles, Smooth Loss, Continuous Treatment, Non-Linear Effect} \\
    \cmidrule(r){2-8}
    & $R^2(\gamma)$ & $R^2(\alpha)$ & MAE & bias & sd & se/sd & covg. \\
    \midrule
    \multicolumn{8}{l}{$n = 2,000$} \\
    Main Spec. & 0.976 & 0.926 & 0.036 & 0.013 & 0.043 & 0.933 & 0.926 \\

    \addlinespace
    \multicolumn{8}{l}{$n = 10,000$} \\
    Main Spec. & 0.990 & 0.964 & 0.016 & 0.002 & 0.019 & 0.947 & 0.938 \\
 
    \bottomrule
    \end{tabular}
\end{table}

\clearpage

\begin{appendix}
\section{Proofs}
\begin{proof}[Proof of Theorem \ref{thm1}]
    This is a special case of Theorem \ref{thm1vunknown} with $v_\rho(X) = \hat{v}_\rho(X) = -1$ and $\epsilon_{\rho n} = 0$. We prove the more general version below.
\end{proof}

\begin{proof}[Proof of Corollary \ref{cor1}]
An upper bound for the critical radius of a MLP
neural net is given in equation (A.10) of \cite{farrell2021deep}. Using the fact that the number of parameters given there is bounded by
$CK^{2}m$ it follows that
\begin{equation}
\delta_{n}\leq C\sqrt{\frac{K^{2}m^{2}\ln(K^{2}m)\ln(n)}{n}},\label{crit rad}%
\end{equation}
where $C$ denotes a generic positive constant. Let $\epsilon_{n}=\inf
_{\alpha\in\mathcal{A},x\in\mathcal{X}}\left\vert \alpha(x)-\alpha
_{0}(x)\right\vert .$ It follows by the uniform approximating bounds given in
FLM, in particular in the first inequality on the top of p. 206, that%
\[
K^{2}m^{2}\ln(K^{2}m)\leq C\epsilon_{n}^{-2d/\beta}(\ln(1/\epsilon_{n}%
)+1)^{7}.
\]
It follows that for any $\varepsilon>0$ and $n$ large enough,
\[
\epsilon_{n}\leq C\{Km\sqrt{\ln(K^{2}m)}\}^{-\beta/d+\varepsilon},
\]
where the presence of $\varepsilon$ allows us to ignore the $(\ln
(1/\epsilon_{n})+1)^{7}$ term. It follows that%
\begin{equation}
\left\Vert \alpha^{\ast}-\alpha_{0}\right\Vert \leq C\epsilon_{n}\leq
C\{Km\sqrt{\ln(K^{2}m)}\}^{-\beta/d+\varepsilon}.\label{approx rate}%
\end{equation}
The conclusion then follows from Theorem 1 and squaring and plugging in the
inequalities from equations (\ref{crit rad}) and (\ref{approx rate}).
\end{proof}

\begin{proof}[Proof of Theorem \ref{thm2}]
    This is a special case of Theorem \ref{thm2vunknown} with $\rho(W, \gamma) = Y - \gamma(X)$. Note that Assumption \ref{ass3} implies Assumption \ref{ass5} (i). Assumptions \ref{ass5} (ii), (iii) and \ref{ass6} are obviously satisfied for this choice of $\rho(W, \gamma) = Y - \gamma(X)$.
\end{proof}

\begin{proof}[Proof of Corollary \ref{cor2}]
By Corollary \ref{cor1}, $\left\Vert \hat{\alpha}_\ell - \alpha_0 \right\Vert = O_p(\epsilon_{\alpha n})$, which satisfies the rate conditions of Assumption \ref{ass4}. The conclusion follows by Theorem \ref{thm2}.
\end{proof}

\begin{proof}[Proof of Theorem \ref{thm1vunknown}] Throughout this proof, let $C > 0$ denote a generic constant (possibly different each time it appears), and let ${\mathbb{E}}_{n}[\cdot]$ denote the
empirical expectation over a sample of size $n$, i.e. ${\mathbb{E}}%
_{n}[Z]=\frac{1}{n}\sum_{i=1}^{n}Z_{i}$,%
\vsedit{\begin{align*}
L_{n}(\alpha, v) &  =~\mathbb{E}_{n}[-2m(W,\alpha)-v(W)\alpha(X)^{2}],\\
L(\alpha, v) &  =~\E[-2m(W,\alpha)-v(W)\alpha
(X)^{2}]
\end{align*}}
Note that
\[
\hat{\alpha}=\operatorname{arg}\min_{\alpha\in\mathcal{A}_n}%
L_{n}(\alpha, \hat{v}_{\rho}).
\]
Since $\Gamma$ is a closed linear space and $\alpha_0$ is defined as the minimizer of $L(a)$ over $\Gamma$, 
then we have the first-order condition that for all $\nu \in \Gamma$:
\begin{align}
    \frac{\partial}{\partial \tau} L(\alpha_0 + \tau\nu, v_\rho)\bigg|_{\tau=0} = 0
\end{align}
Moreover, note that by the linearity of the moment $m$ and linearity of expectation:
\vsedit{\begin{align}
    \frac{\partial}{\partial \tau} L(\alpha_0 + \tau\nu, v)\bigg|_{\tau=0} =~& \E\left[-2m(W;\nu) - 2 v(W)\alpha_0(X) \nu(X)\right]\nonumber\\
    =~& \E\left[-2m(W;\nu) - 2 \E[v(W)\mid X]\alpha_0(X) \nu(X)\right]
\end{align}}
Thus we have:
\vsedit{\begin{align*}
     \frac{\partial}{\partial \tau} L(\alpha_0 + \tau\nu, v)\bigg|_{\tau=0} =~& \frac{\partial}{\partial \tau} L(\alpha_0 + \tau\nu, v)\bigg|_{\tau=0} - \frac{\partial}{\partial \tau} L(\alpha_0 + \tau\nu, v_\rho)\bigg|_{\tau=0}\\
     =~& -2 \E[\alpha_0(X)\, (\E[v(W)\mid X]- v_\rho(X))\, \nu(X)]
\end{align*}}
\vsedit{Define:
$$\hat{\hat{v}}_\rho(X) = \E[\hat{v}_\rho(W)\mid X]$$}
By a Taylor expansion, with $\nu=\alpha-\alpha_0$ and some $\bar{\tau}\in [0,1]$:
\vsedit{\begin{align*}
    L(\alpha, \hat{v}_\rho) - L(\alpha_0, \hat{v}_\rho) =~& \frac{\partial}{\partial \tau} L(\alpha_0 + \tau\nu, \hat{v}_\rho)\bigg|_{\tau=0} + \frac{\partial^2}{\partial \tau^2} L(\alpha_0 + \tau\nu, \hat{v}_\rho)\bigg|_{\tau=\bar{\tau}}\\
    =~& -2 \E[\alpha_0(X) (\hat{\hat{v}}_\rho(X) - v_\rho(X))\,\nu(X)] - 2 \E[\hat{\hat{v}}_\rho(X) \nu(X)^2]
\end{align*}}
By a Cauchy-Schwarz inequality and an AM-GM inequality and since $\alpha_0(X)$ is bounded:
\vsedit{\begin{align*}
    |2 \E[\alpha_0(X) (\hat{\hat{v}}_\rho(X) - v_\rho(X))\,\nu(X)]|\leq~& C\, \|\hat{\hat{v}}_\rho - v_\rho\|\, \|\nu\| \leq \frac{C^2}{2\lambda} \|\hat{\hat{v}}_\rho - v_\rho\|^2 + \frac{\lambda}{2} \|\nu\|^2
\end{align*}}
Moreover, by our assumption $|\hat{v}_\rho(W)| \leq C\implies |\hat{\hat{v}}_\rho(X)|\leq C$. Thus:
\begin{align*}
    2C \|\nu\|^2 \geq - 2 \vsedit{\E[\hat{\hat{v}}_\rho(X) \nu(X)^2]} \geq~& -2\E[v_\rho(X) \nu(X)^2] - \vsedit{|\E[(\hat{\hat{v}}_\rho(X) - v_\rho(X)) \nu(X)^2]|}\\
    \geq~& -2\E[v_\rho(X) \nu(X)^2] - \vsedit{\|\hat{\hat{v}}_\rho - v_\rho\|} \sqrt{\E[\nu(X)^4]}\\
    \geq~& -2\E[v_\rho(X) \nu(X)^2] - C \vsedit{\|\hat{\hat{v}}_\rho - v_\rho\|} \sqrt{\E[\nu(X)^2]}\\
    =~& -2\E[v_\rho(X) \nu(X)^2] - C \vsedit{\|\hat{\hat{v}}_\rho - v_\rho\|}\, \|\nu\|
\end{align*}
Let $\nu_{\ast} = \alpha - \alpha_{\ast}$ and $\nu_0 = \alpha_{\ast} - \alpha_0$, such that $\nu=\nu_{\ast} + \nu_0$. Since $-v_\rho(X)\geq 0$, and $-\E[v_\rho(X) \nu_{\ast}(X)^2] \geq \lambda \E[\nu_{\ast}(X)^2]$, we have:
\begin{align*}
    -\E[v_\rho(X) \nu(X)^2] =~& \E[|v_\rho(X)|\, \nu_{\ast}(X)^2] 
    + 2\E[|v_\rho(X)|\, \nu_{\ast}(X) \nu_0(X)] + \E[|v_\rho(X)|\,\nu_0(X)^2]\\
    \geq~& \E[|v_\rho(X)|\, \nu_{\ast}(X)^2] 
    - 2\E[|v_\rho(X)|\, |\nu_{\ast}(X) \nu_0(X)|]\\
    \geq~& \E[|v_\rho(X)|\, \nu_{\ast}(X)^2] 
    - \frac{1}{2}\E[|v_\rho(X)| \nu_{\ast}(X)^2] - 2 \E[|v_\rho(X)| \nu_0(X)^2]\\
    \geq~& \E[|v_\rho(X)|\, \nu_{\ast}(X)^2] 
    - \frac{1}{2}\E[|v_\rho(X)| \nu_{\ast}(X)^2] - 2 C\E[\nu_0(X)^2]\\
    \geq~& \frac{1}{2}\E[|v_\rho(X)|\, \nu_{\ast}(X)^2] 
     - 2 C\E[\nu_0(X)^2]\\
     \geq~& \frac{\lambda}{2}\E[\nu_{\ast}(X)^2] 
     - 2 C\E[\nu_0(X)^2]
\end{align*}
Combining the last two inequalities:
\begin{align*}
    2C \|\nu\|^2 \geq - 2 \vsedit{\E[\hat{\hat{v}}_\rho(X) \nu(X)^2]} \geq~& \lambda \|\nu_*\|^2 - 4 C\|\nu_0\|^2 - C\, \vsedit{\|\hat{v}_{\rho} - v_{\rho}\|}\, \|\nu\|\\
    \geq~& \lambda \|\nu_*\|^2 - 4 C\|\nu_0\|^2 - C\, \vsedit{\|\hat{\hat{v}}_{\rho} - v_{\rho}\|}\, (\|\nu_*\| + \|\nu_0\|)\\
    \geq~& \frac{\lambda}{2} \|\nu_*\|^2 - 5C\|\nu_0\|^2 - \left(C + \frac{C^2}{2\lambda}\right) \vsedit{\|\hat{\hat{v}}_\rho - v_\rho\|^2} 
\end{align*}
We conclude that for some constant $C$, for any $\alpha\in \Gamma$:
\begin{align}
    \frac{C}{\lambda} \vsedit{\|\hat{\hat{v}}_\rho - v_\rho\|^2}  + C \|\nu\|^2 \geq  L(\alpha, \hat{v}_\rho) - L(\alpha_0, \hat{v}_\rho) \geq \frac{\lambda}{2} \|\nu_{\ast}\|^2 - C \|\nu_0\|^2 -\frac{C}{\lambda}\, \vsedit{\|\hat{\hat{v}}_\rho - v_\rho\|^2}  \label{eqn:quadratic-loss}
\end{align}
Next, by Lemma~11 of \cite{foster2019orthogonal}, the fact that
\vsedit{$-2m(W,\alpha)-v(W)\alpha(X)^{2}$} is Lipschitz with
respect to the vector $(m(W,\alpha), \vsedit{\sqrt{|v(W)|}} \alpha(X))$ and by choosing $\delta:=\delta_{n}+c_{0}\sqrt{\ln(c_{1}/\zeta)/n}$, where
$\delta_{n}$ is an upper bound on the critical radius of ${\text{star}%
}(\sqrt{\mathcal{V}_n}\cdot (\mathcal{A}-\alpha_{*}))$ and ${\text{star}}(m\circ\mathcal{A}%
-m\circ\alpha_{*})$, then with probability $1-\zeta$, for all $\alpha\in \C{A}_n$ and $v\in \mathcal{V}_n$:
\begin{multline*}
\left\vert L_{n}(\alpha, v)-L_{n}(\alpha_{*}, v)-(L(\alpha, v)-L(\alpha%
_{*},v))\right\vert \\
\leq~  O\left(  \delta\left(  \vsedit{\sqrt{\E[|v(W)| (\alpha(X) - \alpha_{\ast}(X))^2]}}+\sqrt{\E[(m(W,\mathcal{\alpha
})-m(W,\alpha_{*}))^{2}]}\right)  +\delta^{2}\right)
\end{multline*}
By MSE-continuity of the moment and the fact that \vsedit{$|\hat{v}_\rho(W)|$} is upper bounded by a constant:
\begin{align*}
\left\vert L_{n}(\alpha, \hat{v}_{\rho})-L_{n}(\alpha_{*}, \hat{v}_{\rho})-(L(\alpha, \hat{v}_{\rho})-L(\alpha%
_{*},\hat{v}_{\rho}))\right\vert=~  O\left(  \delta\,\sqrt{M}\,\Vert\alpha-\alpha_{*}%
\Vert+\delta^{2}\right)  =: \epsilon_{1}(\alpha)
\end{align*}
Finally, since $\hat{\alpha}=\operatorname{arg}\min_{\alpha%
\in\mathcal{A}_n}\hat{L}_{n}(\alpha)$, we have that:
\begin{equation*}
L_{n}(\hat{\alpha}, \hat{v}_\rho)-L_{n}(\alpha_{\ast}, \hat{v}_\rho)\leq 0 
\end{equation*}
Combined with the concentration inequality, yields:
\begin{align*}
    L(\hat{\alpha}, \hat{v}_\rho)-L(\alpha_{\ast}, \hat{v}_\rho) \leq~& L(\hat{\alpha}, \hat{v}_\rho)-L(\alpha_{\ast}, \hat{v}_\rho) - (L_{n}(\hat{\alpha}, \hat{v}_\rho)-L_{n}(\alpha_{\ast}, \hat{v}_\rho))\leq 
    \epsilon_{1}(\hat{\alpha})
\end{align*}
Invoking Equation~\eqref{eqn:quadratic-loss} at $\alpha=\hat{\alpha}$:
\begin{align}
    \frac{\lambda}{2} \|\hat{\alpha} - \alpha_{\ast}\|^2 \leq~& L(\hat{\alpha}, \hat{v}_\rho) - L(\alpha_0, \hat{v}_\rho) + C\|\nu_0\|^2 + \frac{C}{\lambda} \vsedit{\|\hat{\hat{v}}_\rho - v_\rho\|^2} \nonumber\\
    \leq~& L(\hat{\alpha}, \hat{v}_\rho) - L(\alpha_{\ast}, \hat{v}_\rho) + L(\alpha_{\ast},\hat{v}_\rho) - L(\alpha_0, \hat{v}_\rho) + C\|\nu_0\|^2 + \frac{C}{\lambda} \vsedit{\|\hat{\hat{v}}_\rho - v_\rho\|^2} \nonumber\\ 
    \leq~& L(\hat{\alpha}, \hat{v}_\rho) - L(\alpha_{\ast}, \hat{v}_\rho) + 2C\|\nu_0\|^2 + \frac{2C}{\lambda} \vsedit{\|\hat{\hat{v}}_\rho - v_\rho\|^2}\tag{by Equation~\eqref{eqn:quadratic-loss} at $\alpha=\alpha_{\ast}$}\\ 
    \leq~& \epsilon_1(\hat{\alpha}) + 2C\|\alpha_{\ast}-\alpha_0\|^2 + \frac{2C}{\lambda} \vsedit{\|\hat{\hat{v}}_\rho - v_\rho\|^2} \nonumber
\end{align}
By the AM-GM inequality:
\begin{align*}
    \frac{\lambda}{2} \|\hat{\alpha} - \alpha_{\ast}\|^2 \leq \frac{\lambda}{4} \|\hat{\alpha}-\alpha_{\ast}\|^2 + O\left(\frac{M}{\lambda} \delta^2 + \|\alpha_{\ast}-\alpha_0\|^2 + \frac{1}{\lambda} \vsedit{\|\hat{\hat{v}}_\rho - v_\rho\|^2}\right)
\end{align*}
Re-arranging yields:
\begin{align*}
    \|\hat{\alpha} - \alpha_{\ast}\|^2 \leq  O\left(\frac{M}{\lambda^2} \delta^2 + \frac{1}{\lambda}\|\alpha_{\ast}-\alpha_0\|^2 + \frac{1}{\lambda^2} \vsedit{\|\hat{\hat{v}}_\rho - v_\rho\|^2}\right)
\end{align*}
Finally, note that:
\begin{align*}
    \|\hat{\alpha}-\alpha_0\|^2 \leq~& 2\|\hat{\alpha}-\alpha_{\ast}\|^2 + 2 \|\alpha_{\ast}-\alpha_0\|^2\\
    \leq~& O\left(\frac{M}{\lambda^2} \delta^2 + \left(1 + \frac{1}{\lambda}\right)\|\alpha_{\ast}-\alpha_0\|^2 + \frac{1}{\lambda^2} \vsedit{\|\hat{\hat{v}}_\rho - v_\rho\|^2}\right). \qedhere
\end{align*}
\vsedit{Finally, note that by definition $\|\hat{\hat{v}}_\rho - v_\rho\|=\|\hat{v}_\rho - v_\rho\|_X$.}
\end{proof}

\begin{proof}[Proof of Theorem \ref{thm2vunknown}] Throughout this proof, let $C > 0$ denote a generic constant (possibly different each time it appears). To show the first conclusion we verify
Assumptions 1--3 of \cite[CEINR]{chernozhukov2022locally}, with $g(w,\gamma,\theta)$ and $\phi(w,\gamma,\alpha,\theta)$
there given by $m(w,\gamma)-\theta$ and $\alpha(x)\rho(w,\gamma)$
respectively. By Assumption \ref{ass1} and $\left\Vert \hat{\gamma}_{\ell}-\gamma
_{0}\right\Vert \overset{p}{\to} 0,$
\begin{align}
\int\left\Vert g(w,\hat{\gamma}_{\ell},\theta_{0})-g(w,\gamma_{0},\theta
_{0})\right\Vert ^{2}F_{0}(dw)& =\int\{m(w,\hat{\gamma}_{\ell})-m(w,\gamma
_{0})\}^{2}F_{0}(dw) \notag \\ & \leq M \left\Vert \hat{\gamma}_{\ell}-\gamma_{0}\right\Vert
^{2}\overset{p}{\to}0. \label{Ass1i}
\end{align}

By Assumption \ref{ass5} (i), (iii) and $\left\Vert \hat{\gamma}_{\ell}-\gamma
_{0}\right\Vert \overset{p}{\to} 0$,
\begin{align}
\int\left\Vert \phi(w,\hat{\gamma}_{\ell},\alpha_{0},\theta_{0}%
)-\phi(w,\gamma_{0},\alpha_{0},\theta_{0})\right\Vert ^{2}F_{0}%
(dw) &= \int \alpha_{0}(x)^{2}\{ \rho(w, \hat{\gamma}_{\ell})-\rho(w,\gamma_{0})\}^{2}F_{0}(dw) \notag \\
&  \leq C \int \{ \rho(w, \hat{\gamma}_{\ell})-\rho(w,\gamma_{0})\}^{2}F_{0}(dw)
\overset{p}{\to}0. \label{Ass1ii}
\end{align}

By Assumption \ref{ass5} (i), since $\left\Vert \hat{\alpha
}_{\ell}-\alpha_{0}\right\Vert \overset{p}{\to}0$, iterated
expectations gives%
\begin{align}
\int\left\Vert \phi(w,\gamma_{0},\hat{\alpha}_{\ell},\tilde{\theta}_{\ell
})-\phi(w,\gamma_{0},\alpha_{0},\theta_{0})\right\Vert ^{2}F_{0}%
(dw) 
&  =\int\{\hat{\alpha}_{\ell}(x)-\alpha_{0}(x)\}^{2}\rho(W,\gamma_0)^{2}%
F_{0}(dw) \notag \\ & \leq C\left\Vert \hat{\alpha}_{\ell}-\alpha_{0}\right\Vert
^{2}\overset{p}{\to}0. \label{Ass1iii}
\end{align}
Therefore, Assumption 1 (i), (ii), and (iii) of CEINR is satisfied.

Next note that:
\begin{align*}
\hat{\Delta}_{\ell}(w)  &  :=\phi(w,\hat{\gamma}_{\ell},\hat{\alpha}_{\ell
},\tilde{\theta}_{\ell})-\phi(w,\gamma_{0},\hat{\alpha}_{\ell},\tilde{\theta
}_{\ell})-\phi(w,\hat{\gamma}_{\ell},\alpha_{0},\theta_{0})+\phi(w,\gamma
_{0},\alpha_{0},\theta_{0})\\
&  = \{\hat{\alpha}_{\ell}(x)-\alpha_{0}(x)\}\{\rho(w, \hat{\gamma}_{\ell})-\rho(w,\gamma_{0})\}. 
\end{align*}

Let $\bar{\rho}(X,\gamma)=\E[\rho(W,\gamma)\mid X].$ Then by iterated expectations,
the Cauchy-Schwartz inequality, and Assumptions \ref{ass5} and \ref{ass4},
\begin{align}
\int\hat{\Delta}_{\ell}(w)F_{0}(dw)  &  =\int\{\hat{\alpha}_{\ell}(x)-\alpha_{0}(x)\}\{\rho(w, \hat{\gamma}_{\ell})-\rho(w,\gamma_{0})\} F_{0}(dx) \notag \\
& \leq \left\Vert \hat{\alpha}_{\ell}-\alpha_{0}\right\Vert \left\Vert
\bar\rho(\cdot, \hat{\gamma}_{\ell})-\bar\rho(\cdot,\gamma_{0})\right\Vert \notag \\
&  \leq C \left\Vert \hat{\alpha}_{\ell}-\alpha_{0}\right\Vert \left\Vert
\hat{\gamma}_{\ell}-\gamma_{0}\right\Vert = o_{p}(n^{-1/2}). \label{Ass1iv}
\end{align}

Since $\hat{\alpha}_\ell(x)$ and $\alpha_{0}(x)$ are bounded,%
\begin{align}
\int\left\Vert \hat{\Delta}_{\ell}(w)\right\Vert ^{2}F_{0}(dw)  &  =\int%
\{\hat{\alpha}_{\ell}(x)-\alpha_{0}(x)\}^2\{\rho(w, \hat{\gamma}_{\ell})-\rho(w,\gamma_{0})\}^2F_{0}(dw) \notag \\
&  \leq C \left\Vert
\hat{\gamma}_{\ell}-\gamma_{0}\right\Vert^{2} \overset{p}{\to}0, \label{Ass1v}
\end{align}
as in equation \eqref{Ass1ii}. By equations \eqref{Ass1iv} and \eqref{Ass1v} it follows that Assumption 2 (i) of CEINR is satisfied.

Assumption 3 of CEINR follows by Assumption \ref{ass6}. Therefore each of
Assumptions 1--3 of CEINR are satisfied, so the first conclusion follows by
Lemma 15 of CEINR and the Lindeberg-Lévy central limit theorem.

Finally, by the first conclusion $\hat{\theta}\overset{p}{\to}\theta_{0}$ and thus \[\int\{m(w,\hat{\gamma}_{\ell})-\hat{\theta}%
-m(w,\gamma_{0})+\theta_0\}^{2}F_{0}(dw)\overset{p}{\to}0,\] so that the hypotheses of Lemma 16 of CEINR are satisfied, giving the second conclusion. 
\end{proof}

\begin{proof}[Proof of \ref{thm1nonlin}]
he proof would be identical to Theorem~\ref{thm1vunknown}, with the only difference being that $v$ now contains two nuisance parameters $(D, v_\rho)$ and:
\begin{align*}
     \frac{\partial}{\partial \tau} L(\alpha_0 + \tau\nu, \hat{v})\bigg|_{\tau=0} =~& \frac{\partial}{\partial \tau} L(\alpha_0 + \tau\nu, \hat{v})\bigg|_{\tau=0} - \frac{\partial}{\partial \tau} L(\alpha_0 + \tau\nu, v_0)\bigg|_{\tau=0}\\
     =~& \E[\hat{D}(W;\nu) - D(W;\nu)]-2 \E[\alpha_0(X)\, (\vsedit{\E[v(W)\mid X]}- v_\rho(X))\, \nu(X)]
\end{align*}
The first part can then be bounded as:
\begin{align*}
    |\E[\hat{D}(W;\nu) - D(W;\nu)]| \leq \epsilon_{mn} \|\nu\|
\end{align*}
The proof then follows identically to the proof of Theorem~\ref{thm1vunknown}.
\end{proof}

\begin{proof}[Proof of \ref{thm2nonlin}]
It follows exactly as in the proof of Lemma 15 of
CEINR that for each $j$
\begin{align*}
&  \frac{1}{\sqrt{n}}\sum_{i\in I_{\ell}}[\hat{\alpha}_{j\ell}(X_{ji})\rho
_{j}(W_{i},\hat{\gamma}_{j\ell})-\alpha_{j0}(X_{ji})\rho_{j}(W_{i},\gamma
_{j0})]\\
&  =\frac{1}{\sqrt{n}}\sum_{i\in I_{\ell}}[\hat{\alpha}_{j\ell}(X_{ji}%
)-\alpha_{j0}(X_{ji})]\rho(W_{i},\gamma_{j0})+\frac{1}{\sqrt{n}}\sum_{i\in
I_{\ell}}\alpha_{j0}(X_{ji})[\rho(W_{i},\hat{\gamma}_{j\ell})-\rho
(W_{i},\gamma_{j0})]+o_{p}(1)\\
&  =\frac{n_{\ell}}{\sqrt{n}}\int\alpha_{j0}(x_{j})[\rho(w,\hat{\gamma}%
_{j\ell})-\rho(w,\gamma_{j0})]F_{0}(dw)+o_{p}(1),
\end{align*}%
\[
\frac{1}{\sqrt{n}}\sum_{i\in I_{\ell}}[m(W_{i},\hat{\gamma}_{\ell}%
)-m(W_{i},\gamma_{0})]=\frac{n_{\ell}}{\sqrt{n}}\int[m(w,\hat{\gamma}_{\ell
})-\theta_{0}]F_{0}(dw)+o_{p}(1).
\]
Also by Assumption \ref{ass8} it is the case that $\left\Vert \hat{\gamma}_{j\ell
}-\gamma_{j0}\right\Vert <\varepsilon$ for all $j$ with probability
approaching one, so that by the triangle inequality and Assumption 9 iii) we
have%
\begin{align*}
&  \left\vert \frac{1}{\sqrt{n}}\sum_{i\in I_{\ell}}[m(W_{i},\hat{\gamma
}_{\ell})-\theta_{0}+\sum_{j=1}^{J}\hat{\alpha}_{j\ell}(X_{ji})\rho_{j}%
(W_{i},\hat{\gamma}_{j\ell})-\psi(W_{i},\gamma_{0},\alpha_{0},\theta
_{0})]\right\vert \\
&  \leq\frac{n_{\ell}}{\sqrt{n}}\left\vert \int[m(w,\hat{\gamma}_{\ell
})-\theta_{0}+\sum_{j=1}^{J}\alpha_{j0}(x_{j})\rho_{j}(w,\hat{\gamma}_{j\ell
})]F_{0}(dw)\right\vert +o_{p}(1)\\
&  \leq\sqrt{n}C\sum_{j=1}^{J}\left\Vert \hat{\gamma}_{j}-\gamma
_{j0}\right\Vert ^{2}=\sqrt{n}o_{p}((n^{-1/4})^{2})=o_{p}(1),
\end{align*}
where
\[
\psi(w,\gamma_{0},\alpha_{0},\theta_{0}):=m(w,\gamma_{0})-\theta_{0}%
+\sum_{j=1}^{J}\alpha_{j0}(x_{j})\rho_{j}(w,\gamma_{j0}).
\]
The first conclusion then follows by the triangle inequality and the central
limit theorem. The second conclusion follows in analogous way, treating each
$j$ separately, using the arguments in Lemma 16 of CEINR.
\end{proof}

\section{Hyperparameters \label{sec:hyperpar}}
Here we give details on the hyperparameters in the architecture and training of neural nets used in our main specification in Sections \ref{sec:empir} and \ref{sec:simul}.

The regression learner $\hat{\gamma}$ and the debiasing function learner $\hat{\alpha}$ are both parametrized as neural nets with two hidden layers and ReLU activation function. The width of the hidden layers, the learning rate and the training L2 penalty are tuned on a grid based on the out-of-sample loss on a test set (30\% of the data). We train the parameters of the neural net using the Adam optimizer of PyTorch, with a batch size of 128. During training, we randomly drop some layers out with a dropout probability of 0.05. We also do early stopping to avoid overfitting, where we end the training process if the loss on a separate validation set (also 30\% of the data) decreases by less than $10^{-5}$ in 5 consecutive rounds.
\end{appendix}

\begin{acks}[Acknowledgments]
The authors would like to thank Jinyong Hahn, Andres Santos, and Pragya Sur for helpful comments.
\end{acks}

\begin{funding}
Financial support was provided by NSF Grants SES 1757140 and 2242447.
\end{funding}

\bibliographystyle{imsart-nameyear}
\bibliography{references}       


\end{document}